\documentclass[a4paper,12pt]{article}
\usepackage[latin1]{inputenc}
\usepackage[T1]{fontenc}
\usepackage[english]{babel}
\usepackage{amssymb,amsbsy,amsmath,float,amsthm}
\usepackage{babel,indentfirst}
\usepackage{graphicx}
\usepackage{geometry,color}
\usepackage{fancyhdr}
\usepackage{tabularx}
\usepackage{subfigure}
\usepackage{setspace}

\newtheorem{theorem}{Theorem}[section]
\newtheorem{proposition}[theorem]{Proposition}

\newtheorem{remark}[theorem]{Remark}

\newcommand{\supp}{ \rm supp\,}
\newcommand{\sech}{\rm sech\,}
\newcommand{\diag}{\rm diag\,}

\addto\captionsenglish{}

\title{On Semi-Classical Questions Related to Signal Analysis}

\author {B. Helffer\thanks{Université Paris Sud,  Département de Mathématiques,
Bâtiment 425, 91405 Orsay Cedex, France}  and T.M. Laleg-Kirati \thanks{Corresponding author.
Address: INRIA Bordeaux sud-Ouest, Université de Pau et des Pays de l'Adour, UFR Sciences, Bât 1, 64013 Pau, France.
Tel: 0033540175153.
Fax: +33540175150
E-mail: Taous-Meriem.Laleg@inria.fr}}

\date{}
\begin{document}

\maketitle

\paragraph{\textbf{Abstract:}} This study explores the reconstruction of a signal using spectral quantities associated with some self-adjoint realization of an h-dependent Schr\"{o}dinger operator $-h^2 \dfrac{d^2}{dx^2} - y(x)$, $h>0$, when the parameter $h$ tends to $0$. Theoretical results in semi-classical analysis are proved. Some numerical results are also presented. We first consider as a toy model the $\sech^2$ function. Then we study a real signal given by arterial blood pressure measurements. This approach seems to be very promising in signal analysis. Indeed it provides new spectral quantities that can give relevant information on some signals as it is the case for arterial blood pressure signal.

\paragraph{\textbf{Keywords:}} Semi-classical analysis, Schr\"{o}dinger operator, signal analysis, arterial blood pressure

\section{Introduction}

Let $y$ be a positive real valued function on a bounded open interval
$\Omega =]a,b[ $ representing the signal to be analyzed. Following the idea in \cite{Laleg:08}, \cite{LaCrSo:10}, in this study we interpret the signal $y$ as a multiplication operator, $\phi \rightarrow y \cdot \phi$, on some function space.  The spectrum of a regularized version of this operator, i.e. the
 Dirichlet  realization in $\Omega$ or the periodic realization   of
\begin{equation}\label{schrchap2}
H_{h,y} = -h^2\frac{d^2}{dx^2} - y,
 \end{equation}
for a small $h$, is then used for the analysis of $y$  (see \cite{Laleg:08}, \cite{LaCrSo:10}). We will denote by $H_{h,y} ^D$ and $H_{h,y} ^{per}$ these realizations and by $H_{h,y}^{\star}$ one of these two realizations.

We are interested in analyzing the potential $y$ on a compact $K
\subset \Omega$ using the negative eigenvalues $ \lambda_{hn}$,
$n=1,\cdots, N_h$ and some associated orthonormal basis of  real eigenfunctions
$\psi_{hn}$, $n\in \mathbb N $ of  the Schr\"odinger operator $H_{h,y}^*$.
For this purpose, we choose  $\lambda <0$ such that $\lambda$ is a noncritical value of $- y$, $\lambda < \inf \left(-y(a), - y(b))\right)$ and  $- y$ is strictly
  less  than $\lambda$ on  $K$. Then, we show that $y(x)$ can be
  reconstructed in $K$ using the following expression
\begin{equation}\label{formule}
y_{h}(x,\lambda)  = -\lambda +  4 h \sum_{\lambda_{hn} < \lambda} ( \lambda-\lambda_{hn})^{\frac{1}{2}}\, \psi_{hn}(x)^2\,, \quad x\in K.
\end{equation}

This expression is different from the one considered in \cite{LaCrSo:10} which is given by:
\begin{equation}\label{formule2}
y_h(x,0) =    4 h
\sum_{\lambda_{hn } <0 } (-\lambda_{hn})^{\frac{1}{2}}\, \psi_{hn} (x)^2\,,
\quad x\in \mathbb R\,,
\end{equation}
that is corresponding to a problem which is defined, at least at the theoretical level,
on the whole line, under the assumption that  $y$ tends sufficiently rapidly to $0$ at $\infty$.\\
Indeed, in \cite{LaCrSo:10}, the whole potential is recovered, using all
the negative eigenvalues and the associated eigenfunctions of the selfadjoint realization of
$H_{h,y}$ on the line, on the basis of a scattering formula due to  Deift-Trubowitz \cite{DeTr:79}. This formula involves a
remainder whose smallness as $h\rightarrow 0$ is unproved. In addition
the numerical computations are actually done for a spectral problem in
a bounded interval, hence it seems more natural to look directly at
such a problem. Another advantage of our new approach is that we can consider in the same way,
 with $y_{h}(x,\lambda)  = y_{h,\frac 12}(x,\lambda) $,
 \begin{equation}\label{formule}
y_{h,\gamma}(x,\lambda)  = -\lambda +  \big(\frac{ h}{L_\gamma^{cl}} \sum_{\lambda_{hn} < \lambda} ( \lambda-\lambda_{hn})^{\gamma} \, \psi_{hn}(x)^2\big) ^{\frac{2}{2\gamma+1}}, \quad x\in K.
\end{equation}
for $\gamma \in[0,+\infty[$ and $L_\gamma^{cl}$ is a suitable universal semi-classical constant
 and discuss the speed of convergence as $h\rightarrow  0$ of  $y_{h,\gamma}(x,\lambda)$  in function of $\gamma$.\\
 At the theoretical level, we can either choose the Dirichlet  or the periodic realization. Agmon's estimates  (\cite{He,HMR})  show indeed that
 \begin{equation}
 y_{h}^D(x,\lambda) - y_{h}^{per}(x,\lambda) = \mathcal O (h^\infty)\,,
 \end{equation}
 and actually is exponentially small as $h\rightarrow 0$. Nevertheless, at the numerical level, the periodic problem seems to give better results in term of accuracy and convergence speed.\\
In addition, if  we are interested in the analysis of the signal in a specific interval,  the introduction of  a suitable $\lambda$ depending on this interval
enables us  to have a good estimate  of this  part of the signal with a smaller
number of negative eigenvalues. This can have interesting applications
in signal analysis where sometimes the interest is focused on the analysis
of a small part of the signal. Let  us also mention that the interest of our approach is not
really in the reconstruction of the signal that we already have in fact but in
computing some new spectral quantities  that provide relevant
information on the signal. These quantities  could be  the negative
eigenvalues or some Riesz means  of these eigenvalues.
The main application in this study is in the arterial blood pressure (ABP) waveform analysis (the signal is then the pressure).  We refer to  \cite{LaCrSo:10}, \cite{LaMePaCoVa:10}, \cite{LaMeCoSo:07} where for instance it is shown for $\lambda=0$ how these quantities permit  to discriminate between different pathological or physiological situations and also to provide some information on  cardiovascular parameters of great interest.

As described in \cite{LaCrSo:10}, the parameter $h$ plays an important
role in this approach. Indeed, as $h$ becomes smaller,  the approximation
of $y(x)$ by $y_h(x,0)$ improves. We have the same remark in this study and we
will prove the pointwise convergence of $y_{h}(x,\lambda)$ (or more generally $y_{h,\gamma}(x,\lambda)$ ) to $y(x)$ when
$h\rightarrow 0$. This explains our terminology ``semi-classical''.  Note finally that in the applications the choice of $h$ is not necessarily very small. Hence the theoretical analysis given in this article has only for object to give some information for the choice of an optimal $h$ and possibly some $\lambda$.

\section{Main results}
Let  us now present our main result which will be later obtained as a particular case of a more general theorem.

\begin{theorem}\label{maintheorem}~\\
Let  $y$ be a real valued $C^\infty$ function on a bounded open set
$\Omega$. Then,  for  any pair $(K,\lambda)$ such that $K$ is compact and
\begin{equation}\label{asslambdaK}\left\{
\begin{array}{l}
\lambda < \inf (-y (a), - y(b)) \,,\\
 y(K) \subset ]- \lambda, +\infty[\,, \\
-\lambda \mbox{ is not a critical value of } y\,,
\end{array}\right.
\end{equation}
and, uniformly  for  $x\in K$, we have
\begin{equation}
y(x) = - \lambda + \lim_{h\rightarrow 0} {4 h \sum_{\lambda_{hn} < \lambda} (\lambda-\lambda_{hn})^{\frac{1}{2}} \, \psi_{hn} (x) ^2} \,,
  \end{equation}
where $\lambda_{hn}$ and $\psi_{hn}$ denote the  eigenvalues and some
associated $L^2$-normalized  real eigenfunctions of  $H_{h,y}^*$ in $\Omega$. Here we recall that
 $H_{h,y}^*$ denotes either $H_{h,y}^D$ or $H_{h,y}^{per}$.
\end{theorem}
\begin{remark}~\\
We do not impose in the  statement of Theorem  \ref{maintheorem}  the positivity of $y$ and the negativity of
$\lambda$.
If we want to recognize  the signal $y$ on a given compact $K$, then any choice of
$\lambda$ such that \eqref{asslambdaK} is satisfied is possible. There are
at this stage two contrary observations. To choose $\lambda$ as small as possible will give the
advantage that we need less eigenvalues and eigenfunctions to
compute, but this is only  true in the limit $h\rightarrow 0$ and this  imposes
to compute more eigenvalues! For a given $h$, take a larger $\lambda$ seems
numerically better.  An explanation could come from another asymptotic analysis
 which will not be done here. Hence our theorem gives just some light on what should be
done numerically and one has to be careful with their interpretation.
\end{remark}

We will obtain the proof by  using a suitable extension of Karadzhov's theorem on the spectral function. One can suspect that the
 convergence could be better due to the fact that $\lambda$ is not
 critical, but it should not affect so much  the computation if $y(K)
 + \lambda $ is sufficiently large.\\
So the basic idea is that
$$
 h \sum_n \chi (\lambda - \lambda_{hn} )\, \Psi_{hn}(x)^2= (2\pi)^{-1} \int \chi (\lambda +y(x)  -\xi^2 ) \,d\xi +
 \mathcal O(h)\,,
$$
uniformly with respect to $x$ in a  compact and that we can have a complete expansion in $h$ if $\supp \chi \subset\subset ]0,+\infty[$.  In what follows we omit the subscript $y$ in $H_{h,y}$. The main point here (in the case of the line) is, following \cite{HR1},  that $\chi(\lambda -H_h)$ can be considered as an $h$-pseudodifferential operator whose $h$-symbol admits an expansion
 \begin{equation}
 q(x,\xi,h) \sim \sum_{j\geq 0} h^j q_j(x,\xi)
 \end{equation}
  where the $q_j(x,\xi)$ have compact support,
  \begin{equation}
  q_0(x,\xi)=\chi (\lambda + y(x) -\xi^2) \end{equation}
  and
  \begin{equation}
  q_1(x,\xi)=0\,.
  \end{equation}

  This means that, for $u\in C_0^\infty(\mathbb R)$, we can write
  \begin{equation}
\left( \chi(\lambda -H_h) \, u \right)(x) = (2\pi h)^{-1} \int e^{ \frac ih \,(x-y)\cdot \xi } q(\frac{x+y}{2},\xi,h) u(y) dy d\xi\,.
  \end{equation}

 We are then considering the restriction to the diagonal of the distribution kernel $K_h$ of  $\chi(\lambda -H_h)$
  \begin{equation}
   K_h(x,y):= (2\pi h)^{-1} \int e^{ \frac ih \,(x-y)\cdot \xi } q(\frac{x+y}{2},\xi,h) d\xi\,,
  \end{equation}
  which at $(x,x)$ becomes
  \begin{equation}
   K_h(x,x):= (2\pi h)^{-1} \int  q(x,\xi,h) d\xi\,,
  \end{equation}
  and admits the expansion in powers of $h$
  \begin{equation}
   K_h(x,x)\sim (2\pi h)^{-1} \sum_{j\geq 0} h^j \int q_j(x,\xi) d\xi \,.
  \end{equation}
  Moreover, when integrating over $x$, we get the trace of $\chi(\lambda -H_h)$ ~
  \begin{equation}
 {\rm Tr \;}  \chi(\lambda -H_h) =  (2\pi h)^{-1} \int q(x,\xi,h) dx d\xi\,.
  \end{equation}

We also note  for future use (see (2.3) in \cite{HeRo:90a})
  \begin{equation}\label{coeffh2}
  q_2(x,\xi) = \frac 14  \chi''(\lambda + y(x) -\xi^2)  y''(x)
  +  \frac{1}{12}  \chi'''(\lambda + y(x) -\xi^2) (2 y''(x)\xi^2 -  y'(x)^2)\,.
  \end{equation}

The hope is that when $\chi$ is replaced by $\lambda _+^\gamma$, with $\gamma >0$,
 we keep a remainder in $\mathcal O (h^{1+\gamma})$ in the right hand
 side.
This reads
\begin{equation}\label{heuristics1}
 h \sum_n (\lambda -\lambda_{hn})_+^{\gamma} \, \Psi_{hn}(x) ^2=
(2\pi)^{-1} \int  (\lambda -\xi^2 + y(x))_+^{\gamma} d\xi +
 \mathcal O(h^{1+\gamma})\,.
\end{equation}
This suggests the change of variable $\xi =\eta (\lambda + y(x))_+^{\frac
  12}$
 and leads to
\begin{equation}\label{heuristics2}
 h \sum_n (\lambda -\lambda_{hn})_+^{\gamma}\, \Psi_{hn}(x) ^2=
(\lambda  + y(x) )_+^{\gamma + \frac 12}
(2\pi)^{-1} c_\gamma +
 \mathcal O(h^{1+\gamma})\,.
\end{equation}
with
\begin{equation}\label{defcgamma}
c_\gamma :=   \int  (1 -\eta^2)_+^{\gamma} d\eta \,.
\end{equation}
The hope is that the remainder will be uniform for $x$ in
 a compact $K$ such that $y(K) \subset ]-\lambda, +\infty[$.

This theorem  is proven in \cite{Karadzhov2} (with some details given
in
 \cite{Karadzhov1}) for $\gamma=0$, see in the next section.  For our application, we need to have at least the case $\gamma=\dfrac{1}{2}$ but we will also suggest  in our conclusion that considering $\gamma$ larger
  could be better.

\section{Former results }
We recall one of the basic results of Karadzhov \cite{Karadzhov1}.  We come back to the more standard notation by writing $V=-y$.
\begin{theorem}\label{theoremKaradzhov}~\\
Let  $A_h = -h^2 \dfrac{d^2}{dx^2} + V $, with a $C^\infty$ potential $V$ on the line
 with $V$ tending to $+\infty$. Let  $e_h$ be the spectral function associated to $A_h$
such that $\displaystyle{e_h(\lambda,x,y)= \sum_{\lambda_j\leq \lambda
  }{\psi_{hj}(x) \, \psi}_{hj}(y)}$. Then,  for any compact $K$  such that $K\subset V^{-1} (]-\infty,\lambda[)$, we have, for $x\in K$
\begin{equation}\label{estika}
e_h(\lambda,x,x)= \pi ^{-1}  \, (\lambda -V (x)  )^{\frac{1}{2}} h^{-1} + \mathcal O(1), \quad h\rightarrow 0,
\end{equation}
uniformly in $K$.
\end{theorem}

The next theorem was established by Helffer-Robert  \cite{HeRo:90a} in connection with the analysis of the Lieb-Thirring conjecture. This will not be enough because we will need a point-wise estimate  and this is an integrated version, but this indicates in which direction we want to go.

\begin{theorem}\label{theoremehelfferrobert}~\\
Let  $V \in C^\infty(\mathbb{R})$, with
 $$
- \infty < \inf V <\liminf_{|x|\rightarrow +\infty} V \,.
$$
 Let  $\lambda\in
]\inf V,\liminf_{|x|\rightarrow +\infty} V[$  and suppose that $\lambda$ is not a critical value for $V$.
We denote by~:
\begin{equation}\label{moyennes de riesz}
S_\gamma(h,\lambda) = \sum_{\lambda_{hn} \leq
\lambda}{\left(\lambda -
\lambda_{hn}\right)^\gamma},\quad\quad \gamma
\geq 0
\end{equation}
the Riesz means of the eigenvalues $\lambda_{hn}$
less than $\lambda$ of $A_h$. Then for $0\leq \gamma \leq 1$, we have:
\begin{equation}\label{theoreme_moyenne de Riesz equation1}
    S_\gamma(h, \lambda)= \frac{1}{h}\left(L_\gamma^{cl}\int_{-\infty}^{+\infty}{\left(\lambda -
V(x)\right) _+ ^{\gamma +\frac{1}{2}} dx} +
    \mathcal O(h^{1+\gamma})\right)
\end{equation}
where $|\cdot |_+$ is the positive part  and $L_\gamma^{cl}$,
known as the classical Weyl  constant, is given by
\begin{equation}\label{theoreme_moyenne de Riesz equation3}
   L_{\gamma}^{cl}=\frac{\Gamma(\gamma+1)}{2\sqrt{\pi}
     \Gamma(\gamma+\frac{3}{2})}.\end{equation}
     \end{theorem}

     Note that $L_\gamma^{cl}= c_\gamma (2\pi)^{-1}$ and we recall that $\Gamma (\dfrac{3}{2})=\dfrac{\sqrt{\pi}} {2}$. Hence in particular $L_0^{cl} =\dfrac {1}{ \pi}$ and $L_\frac 12^{cl} = \dfrac{1}{4}$.

\section{Proof of the main theorem}

\subsection{Pointwise asymptotics for the Riesz means}
The statement of the main theorem will correspond to the case $\gamma=\dfrac{1}{2}$ of the following
\begin{theorem}\label{theoremKaradzhov}~\\
Let   $H_h^{\star} $ be the realization of $H_h = -h^2\dfrac{d^2}{dx^2} + V $, with a $C^\infty$ potential $V=-y$ (We can either consider $H_h^D$ or $H_h^{per}$). Let  $e_h^\gamma$ be defined by
$$\displaystyle{e_h^\gamma (\lambda,x,y)= \sum_{\lambda_{hj}\leq \lambda
  }{(\lambda -\lambda_{hj} )^\gamma \psi_{hj} (x) \psi}_{hj} (y)}\,.
  $$
  For any pair $(K,\lambda)$ satisfying \eqref{asslambdaK} , then  we have, for $x\in K$,
\begin{equation}\label{estika}
e_h^\gamma (\lambda,x,x)= (2\pi)^{-1}  (\lambda -V (x)  )^{\gamma + \frac{1}{2}} c_\gamma h^{-1}  + \mathcal O(h^{\gamma}), \quad h\rightarrow 0,
\end{equation}
uniformly in $K$.
\end{theorem}
We refer to the heuristics starting from \eqref{heuristics1} for understanding the main term in \eqref{estika}.

As in the proof by Helffer-Robert of Theorem \ref{theoremehelfferrobert} \cite{HeRo:90a}, we can distinguish two steps corresponding to the contribution which is close to $\lambda$ and to the contribution which is "far" from $ \lambda$. This will be done by a cut-off in energy.
\subsection{Far from $\lambda$}
 We consider a function $\chi$ with $\supp \chi
\subset\subset ]-\infty,0[$ and consider the expression
\begin{equation}
f_{h,\chi} ^\gamma(\lambda,x,x)  :=  \sum_{j}  \chi(\lambda_{hj}-\lambda) |\lambda - \lambda_{hj}|^\gamma
  \, \psi_{hj} (x) ^2\,,
  \end{equation}
  and we prove that
  \begin{proposition}~\\
  If  $(\lambda,K)$ satisfies the condition
  \begin{equation}\label{hypaubrd}
  \lambda < \inf (V(a),V(b))\,,
  \end{equation}
  and if $K$ is compact in $]a,b[$ then, uniformly on $K$, we have
  \begin{equation}
  f_{h,\chi} ^\gamma(\lambda,x,x)  \sim (2\pi)^{-1} h^{-1} \left(\sum _{\ell \geq 0} \alpha_{\ell,\chi}^\gamma  (x) h^\ell\right).
  \end{equation}
  where
  \begin{equation}
  \alpha_{0,\chi}^\gamma (x) = \int    \chi(\lambda -\xi^2 - V(x) ) |\lambda -V(x) -\xi^2|^\gamma \,d\xi\,.
  \end{equation}
  and the $\alpha_{\ell,\chi}^\gamma $ are $C^\infty$ functions on $]a,b[$.
  \end{proposition}
  \subsection{Close to $\lambda$}
  We now consider
 \begin{equation}
  g_{h,\chi} ^\gamma(\lambda,x,x)  :=  \sum_{\lambda_{hj} \leq \lambda}  (1 -\chi (\lambda_{hj} -\lambda)) |\lambda - \lambda_{hj}|^\gamma \,\psi_{hj}(x)^2 \,.
  \end{equation}

  This quantity involves only the eigenvalues (and corresponding eigenfunctions) which are close to $\lambda$. Of course, we have
  \begin{equation}
  e_{h}^\gamma (\lambda,x,x) = f_{h,\chi} ^\gamma(\lambda,x,x)   +  g_{h,\chi} ^\gamma(\lambda,x,x)  \,,
  \end{equation}
  and for coming back to our main theorem we keep in mind that
  \begin{equation}
   e_{h}^0 (\lambda,x,x) =  e_{h} (\lambda,x,x) \,.
   \end{equation}

  Assumption \eqref{hypaubrd} together with Agmon estimates permits  (in the two cases, Dirichlet or periodic) to reduce to a global problem on $\mathbb R$ where $V$ is replaced by a new potential $\widetilde V$ coinciding with $V$ on $K$. Moreover, we can impose conditions on $\widetilde V$ at $\infty$ permitting to
  use the global calculus of Helffer-Robert calculus (or the semi-classical Weyl calculus)  and to get the result.

 The support of $\chi$ is now chosen so that $\chi=1$ on $]-\infty,-\epsilon]$
  with $\epsilon >0$ small enough. The second proposition is~
\begin{proposition}~\\
   For any pair $(K,\lambda)$ satisfying \eqref{asslambdaK}  and for a sufficiently small $\epsilon>0$, we have uniformly for $x\in K$
  \begin{equation}
  g_{h,\chi} ^\gamma(\lambda,x,x)  = (2\pi)^{-1} h^{-1} \left(\beta_{0,\chi}^\gamma  (x)  +\mathcal O (h ) \right),
  \end{equation}
  where
  \begin{equation}  \beta_{0,\chi}^\gamma (x) = \int    (1-\chi (\lambda -\xi^2 - V(x) )) (\lambda -V(x) -\xi^2)_+^\gamma \,d\xi\,.
  \end{equation}
 \end{proposition}

  Here the proof is  a rather direct consequence of Karadzhov's paper.
   We are indeed analyzing
   $$
    g_{h,\chi} ^\gamma(\lambda,x,x)  = \int _{\mu \leq \lambda} (1- \chi (\lambda -\mu)) \,(\lambda -\mu) ^{\gamma} \,d e_h(\mu,x,x)\,.
   $$

Using an integration by parts, we obtain
    $$
    g_{h,\chi} ^\gamma(\lambda,x,x)  = \int _{\mu \leq \lambda}  \phi (\lambda -\mu) \, e_h(\mu,x,x) d\mu\,,
   $$
   with
   $$
   \phi(t)= \left( (1-\chi (t)) t_+^\gamma\right)'\,.
   $$

   We can then apply the estimate \eqref{estika} we had for $e_h(\mu,x,x)$  which is uniform for $\mu$ close to $\lambda$ and we obtain by integration over $]-\infty,\lambda]$ (and the reverse integration by parts)
     the result modulo $\mathcal O(h)$.

   \begin{theorem}~\\
   For any pair $(K,\lambda)$ satisfying \eqref{asslambdaK}  and sufficiently small $\epsilon>0$, we have uniformly for $x\in K$
  \begin{equation}
  g_{h,\chi} ^\gamma(\lambda,x,x)  = (2\pi)^{-1} h^{-1} \left(\sum _{0 \leq \ell < \gamma +1} \beta_{\ell,\chi}^\gamma  (x) h^\ell  +\mathcal O (h^{1+\gamma} ) \right),
  \end{equation}
  with
  \begin{equation}
  \beta_{0,\chi}^\gamma (x) = \int    (1-\chi (\lambda -\xi^2 - V(x) )) (\lambda -V(x) -\xi^2)_+^\gamma \,d\xi\,,
  \end{equation}
  and
   \begin{equation}\label{conj3}
  e_h^\gamma(\lambda,x,x)  = (2\pi)^{-1} h^{-1} \left(\sum _{0 \leq \ell < \gamma +1} \alpha_{\ell}^\gamma  (x) h^\ell  +\mathcal O (h^{1+\gamma} ) \right)\,.
  \end{equation}
  with
  \begin{equation}
  \alpha_{0}^\gamma (x) = \int   (\lambda -V(x) -\xi^2)_+^\gamma \,d\xi\,.
  \end{equation}
 Here the $\beta_{\ell,\chi}^\gamma $ and $\alpha_{\ell}^\gamma $ are $C^\infty$ functions in a neighborhood of $K$.
   \end{theorem}

   A direct proof should be given for having this remainder. One should follow Karadzhov's proof (in its easy part because we are far from $V(x)=\lambda$) and improve the Tauberian theorem used in this paper, using the ideas of \cite{HeRo:90a} (in a non integrated version). It is also  based on the approximation of $(1-\chi) H_{h} \exp(-i t \frac{H_{h}}{h})$ by a Fourier Integral Operator.  We do not give the details here.

   \paragraph{More information on the coefficients}~\\
   We get from  \cite{HeRo:90a}  that $\alpha_{1}^\gamma  (x)  =0$ . Actually
    only the integrated version of this claim is given (see (0.12) there and have in mind that the subprincipal symbol vanishes) but coming back to the proof of \cite{HR1} gives the statement .
    It is also proven
    by Helffer-Robert \cite{HeRo:90a} that $\alpha_{2}^\gamma $
    is not identically zero.  The computation is easier for $\gamma >2$. We follow what was done in  \cite{ HeRo:90a} (Equation 2.20),  but we can no more integrate in the $x$ variable, so the simplification obtained in this paper by performing an integration by parts in the $x$ variable is not possible.
   We get, for $\gamma >2 $,  (see \eqref{coeffh2})
    \begin{equation}
    \alpha_{2}^\gamma (x) = L_\gamma (\lambda-V(x))^{\gamma -\frac 32} V''(x)\
     + L'_\gamma (\lambda-V(x))^{\gamma -\frac 52} V'(x)^2,
   \end{equation}
    where $L_\gamma$ and $L'_\gamma$  are universal computable constants.
    So when increasing $\gamma$, we suspect that we improve
     the semi-classical approach of $y(x)$ when we increase $\gamma$ from $0$ to $1$ and that
      we do not improve anymore for larger $\gamma$. Fig. \ref{h0001-gamma05-1-2}  confirms this guess
       for $\gamma =2$. The asymptotic behavior of $y(x)-y_{h,\gamma}(x)$ is indeed given
      by $h^2 \delta_2^\gamma (x)$, where $\delta_2^\gamma$ is a $C^\infty$ function in a neighborhood of $K$ directly computable from  $\alpha_{2}^\gamma (x) $  and $V$. In the case $\gamma=\dfrac{1}{2}$, we see on the contrary some oscillation.

\section{Some numerical examples}

To study the "concrete" validity of our formula, we have performed some numerical tests using MATLAB software. We have chosen to use a pseudo-spectral Fourier method instead of a finite differences method  for the discretization of the problem. Indeed, Fourier method gives better results in term of  accuracy and speed convergence. However Fourier method requires periodic boundary conditions  so the numerical tests have been done on the periodic realization $H^{per}_{h,y}$.

We consider a grid of $M$ equidistant points $x_j$, $j=1,\cdots,M$
such that
\begin{equation}
a=x_1<x_2<\cdots<x_{M-1}<x_M=b.
\end{equation}
We denote $y_j$ and $\psi_j$ the values of $y$ and $\psi$ at the grid points
$x_j$, $j=1,\cdots, M$
\begin{equation}
y_j=y(x_j),\quad \psi_j=\psi(x_j),\quad j=1,\cdots, M.
\end{equation}
Therefore, the discretization of $H^{per}_{h,y}$ leads to the following eigenvalue matrix problem
\begin{equation}\label{sch_discrete2}
    \left( -h^2\boldsymbol{D}_2-  \diag\left( \boldsymbol{Y} \right) \right)\underline{\boldsymbol{\psi}}=-\lambda\underline{\boldsymbol{\psi}},
\end{equation}
where $\diag(\boldsymbol{Y})$ is a diagonal matrix whose elements are $y_j$,
$j=1,\cdots,M$ and $\underline{\boldsymbol{\psi}}= \left[\psi_1 \:\:
\psi_2,\:\:\cdots\:\: \psi_{M-1}\:\: \psi_M\right]^{T}$. $D_2$ is
the second order differentiation matrix for a pseudo-spectral Fourier method \cite{LaCrSo:10}, \cite{Tre:00}. Note that the analysis is only relevant if $M$ is large, and $h$ can not be too small in comparison with the distance $\dfrac{b-a}{M-1}$ between two consecutive points. We denote $N_h$ the number of negative eigenvalues of $H_{h,y}$ and $N_{h,\lambda}$ the number of negative eigenvalues less than $\lambda$.

Different values of the parameter $h$, $\lambda$ and $\gamma$ have been considered, sometimes outside the probable domain of validity of the theoretical analysis. We concentrate our  analysis on  two examples. We first look at the $\sech^2$ function which is a regular function and then we  consider the case of arterial blood pressure signal. The latter is given by measured data at the finger level provided by physicians and it is difficult in this case to speak of regularity of the signal.

\subsection{ The $\sech^2$ signal.}
We first study as toy model the example of the $\sech^2$ function defined  on $[0,10]$ by
\begin{equation}
    y(x) =   \sech^2(x-5) \,.
\end{equation}

It is a well studied potential when considered on the whole line with explicitly known negative spectrum for the associated Schr\"odinger operator for some specific values of $h$.

\begin{figure}[htbp]
\begin{center}
\subfigure[$h=0.01$]{\includegraphics[width=7cm]{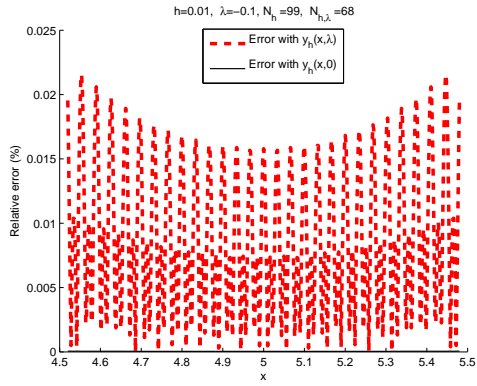}}
\subfigure[$h=0.001$]{\includegraphics[width=7cm]{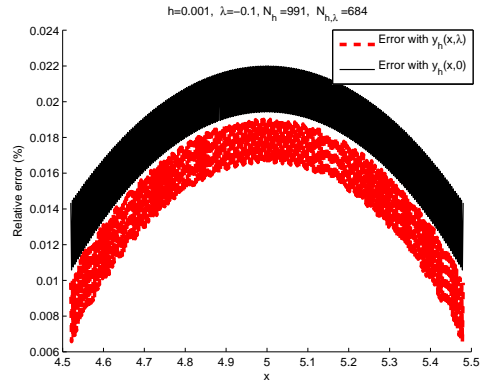}}
\caption{Relative error with  $\lambda=0$ and $\lambda=-0.1$ for the $\sech^2$ example}
\label{erreur-sech-lamda01}
\end{center}
\end{figure}

\begin{figure}[htbp]
\begin{center}
\subfigure[$h=0.01$]{\includegraphics[width=7cm]{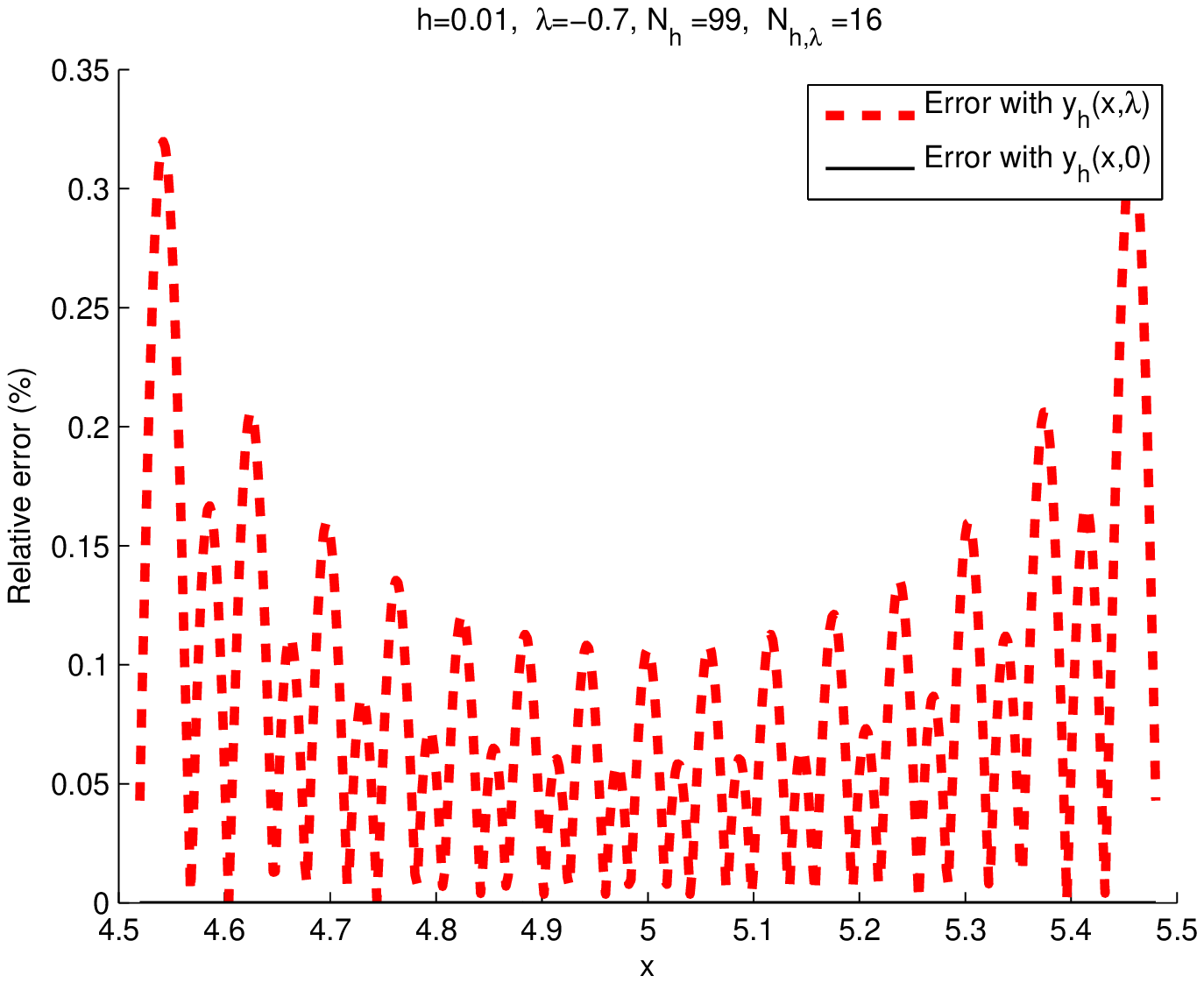}}
\subfigure[$h=0.001$]{\includegraphics[width=7cm]{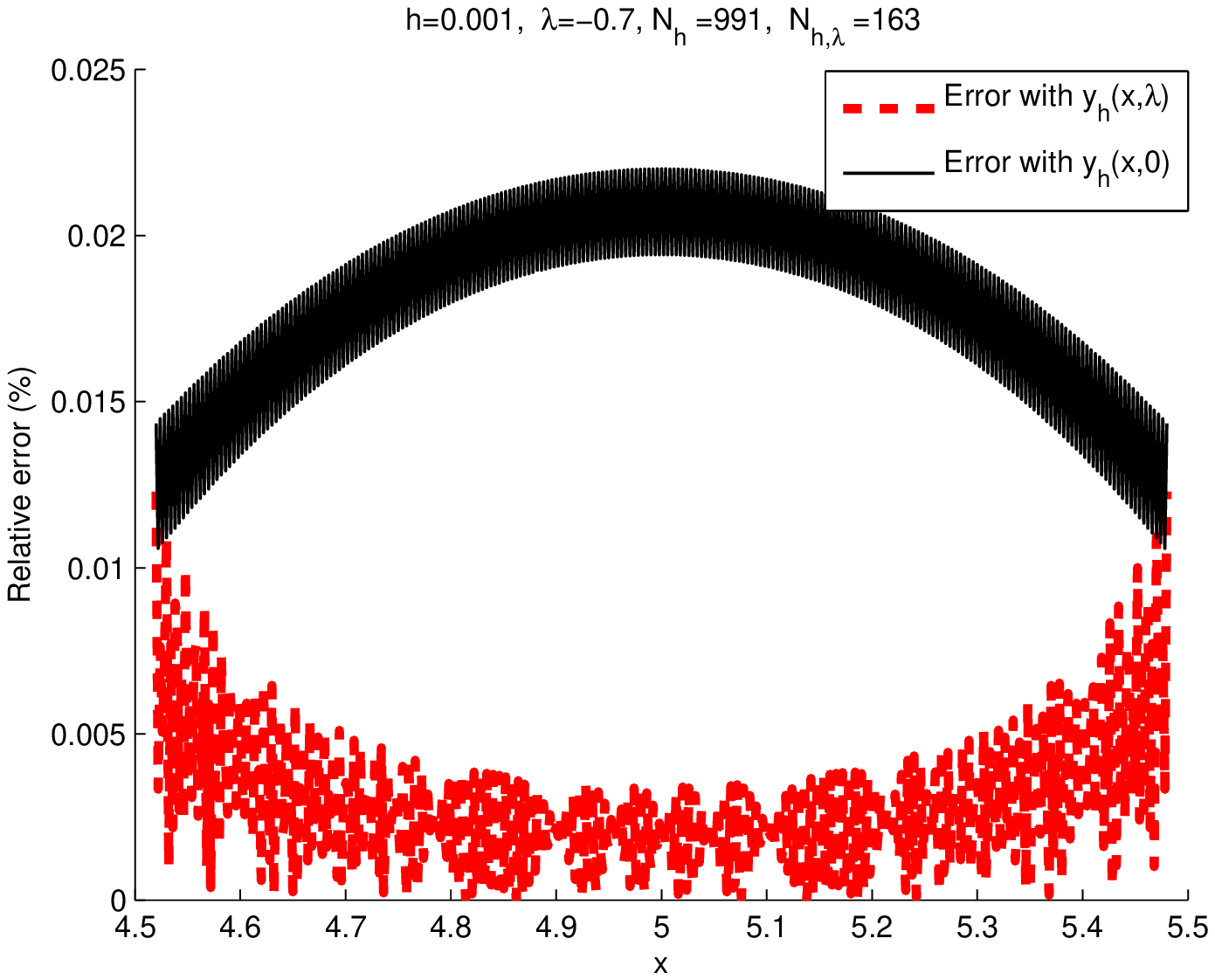}}
\caption{Relative error for $\lambda=0$ and $\lambda=-0.7$ for the $\sech^2$ example}
\label{erreur-sech-lamda07}
\end{center}
\end{figure}

We analyze  the reconstruction of a part of the signal given by $- y <\lambda_1$, with $\lambda_1= -0.8$. For this purpose, we use Formula (\ref{formule2}) and compare the reconstruction of $y(x)$ with this formula with its reconstruction with $y_{h,\gamma}(x,\lambda)$ for different values of $h$ and $\lambda$. Fig. \ref{erreur-sech-lamda01} - fig. \ref{erreur-sech-lamda07} illustrate the results for $\gamma=\dfrac{1}{2}$. We notice that for $h= 0.01$, the reconstruction  of this part of the signal is better with $y_h(x,0) $. However, when $h=0.001$ the reconstruction  is better with $y_{h}(x,\lambda)$.

\begin{figure}[htbp]
\begin{center}
\subfigure[]{\includegraphics[width=7cm]{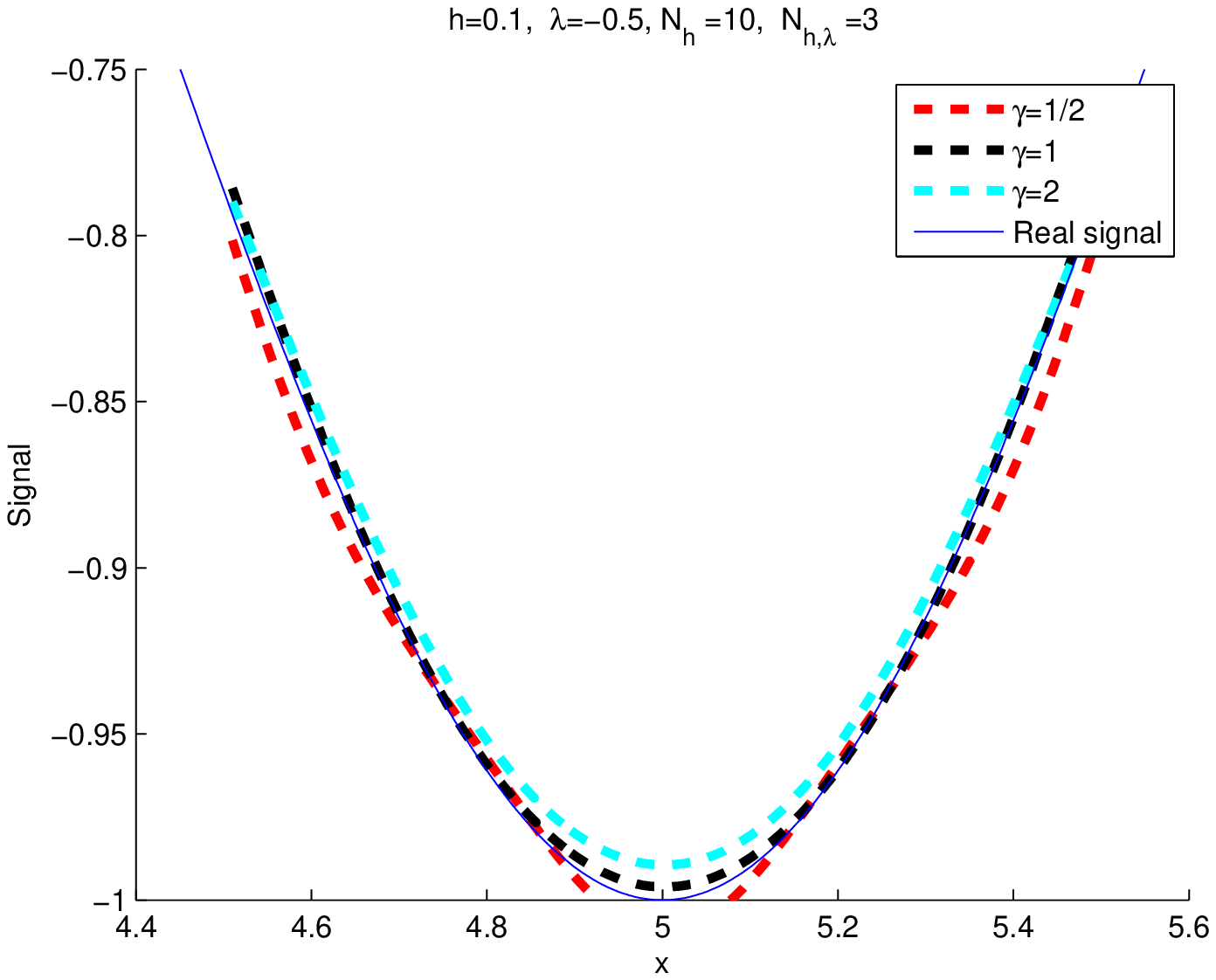}}
\subfigure[]{\includegraphics[width=7cm]{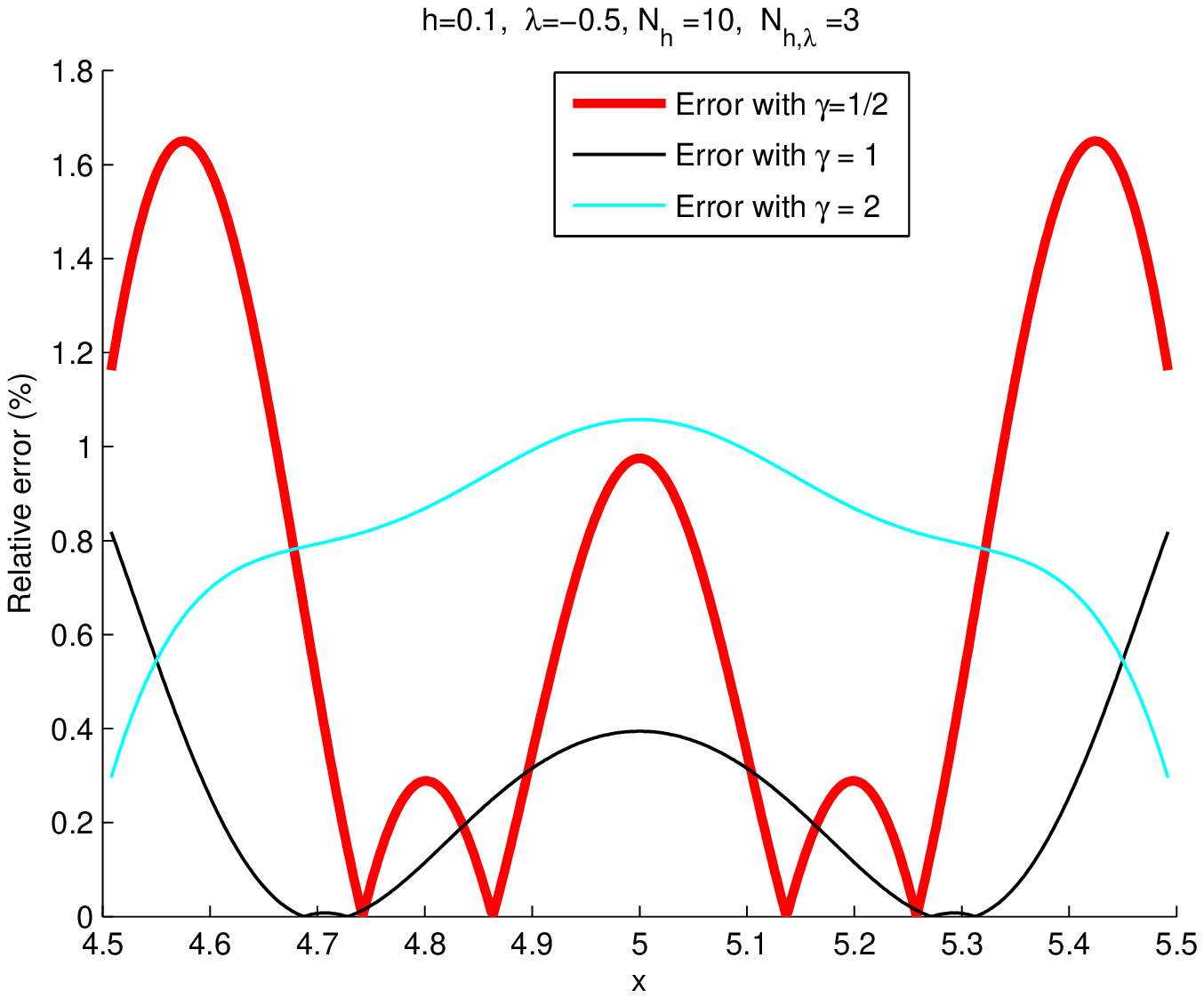}}
\caption{Reconstruction of a part of the $\sech^2$ signal (a) and  relative error (b) for different values of $\gamma$ and $h=0.1$}
\label{h01-gamma05-1-2}
\end{center}
\end{figure}

\begin{figure}[htbp]
\begin{center}
\subfigure[]{\includegraphics[width=7cm]{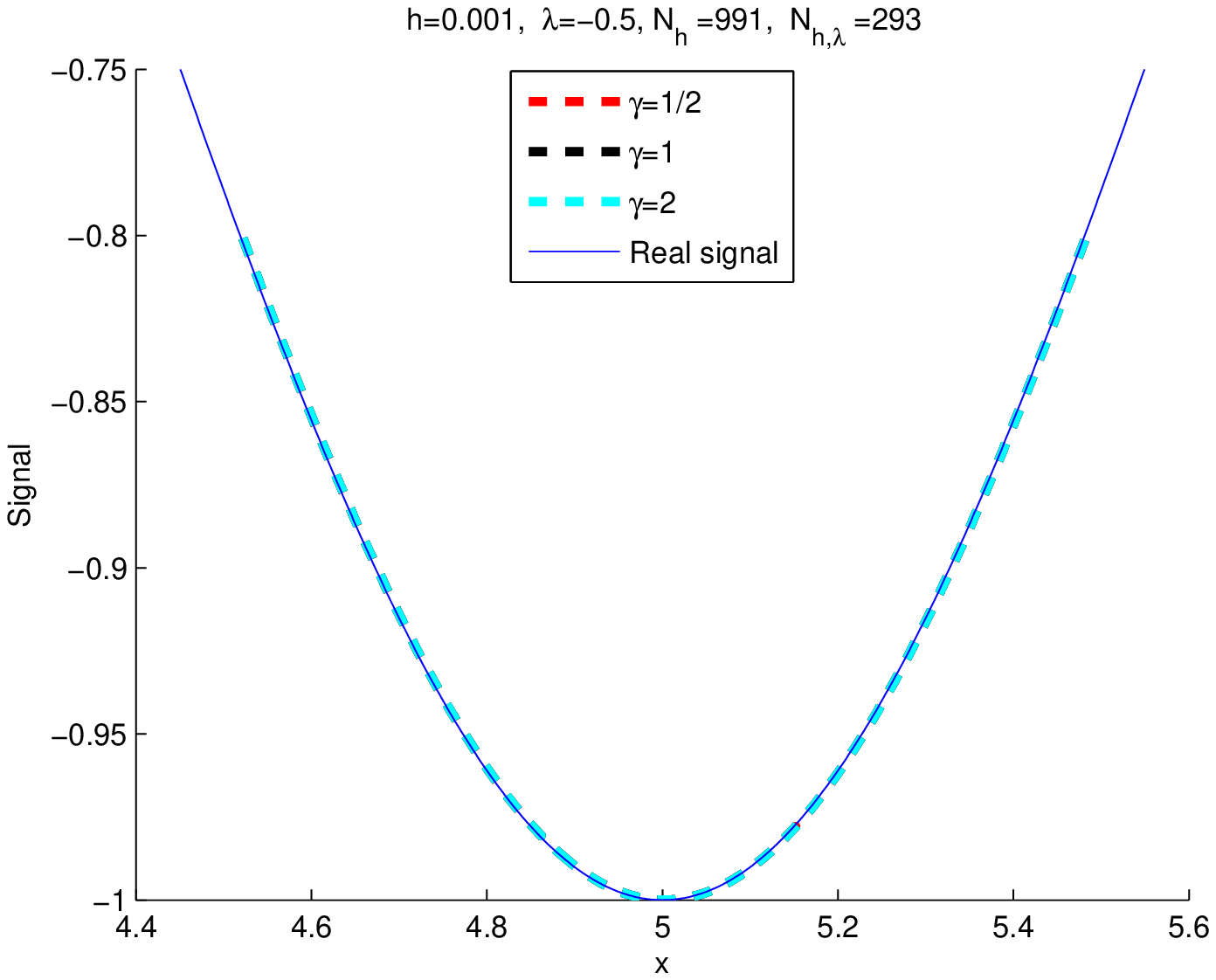}}
\subfigure[]{\includegraphics[width=7cm]{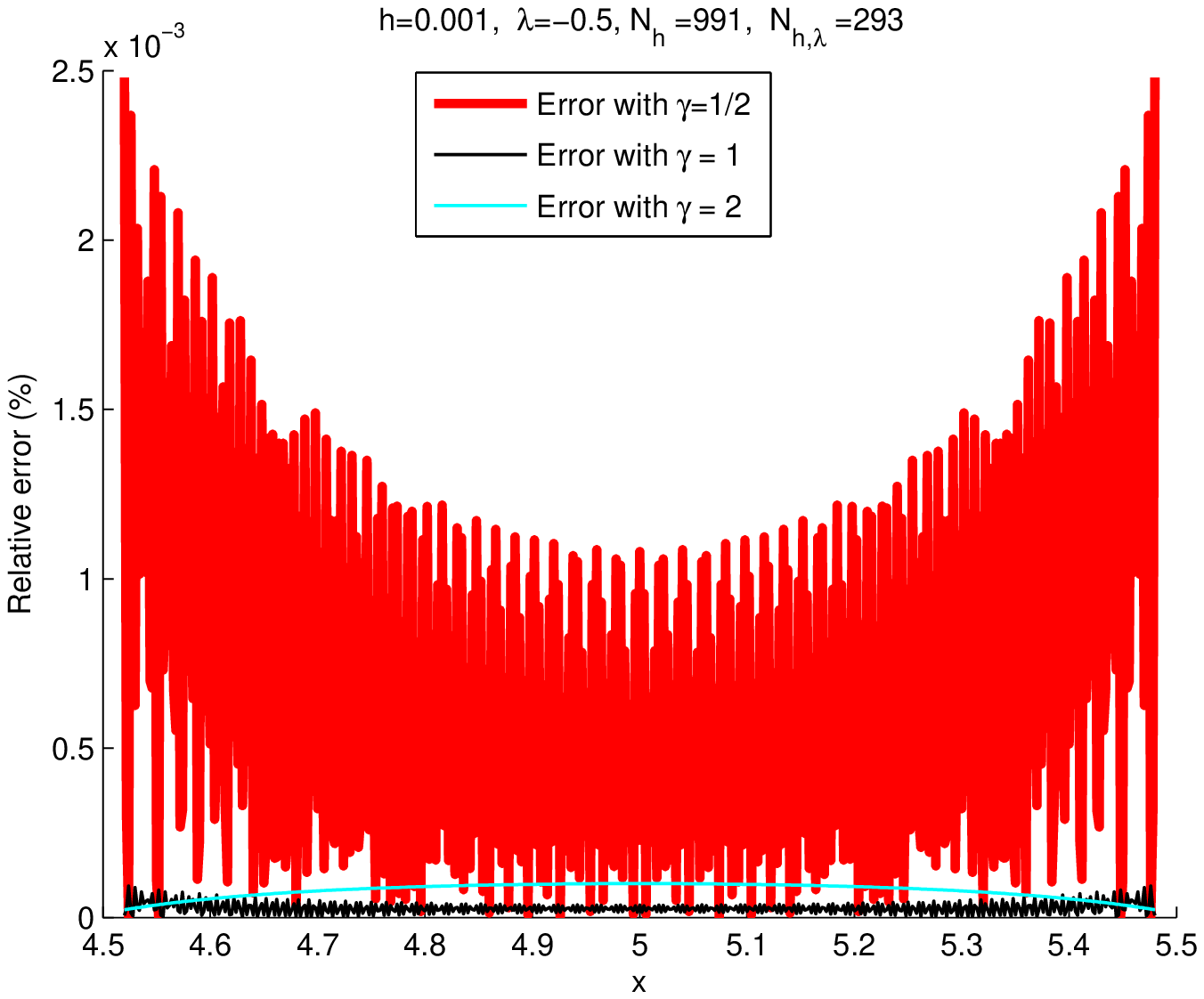}}
\caption{Reconstruction of a part of the $\sech^2$ signal (a) and relative error (b) for different values of $\gamma$ and $h=0.001$}
\label{h0001-gamma05-1-2}
\end{center}
\end{figure}

Then, we analyze the error  with different values of $\gamma$ and $h$, $\lambda$ being fixed. Fig. \ref{h01-gamma05-1-2} - fig. \ref{h0001-gamma05-1-2} illustrate the results for $\lambda=-0.5$. The optimal choice of $\gamma$ for fixed $h$ and $\lambda$  seems to be $\gamma = 1$. Note also that the error for $\gamma=2$ is regular as it was explained previously.

\begin{figure}[htbp]
\begin{center}
\subfigure[$h=0.01$]{\includegraphics[width=7cm]{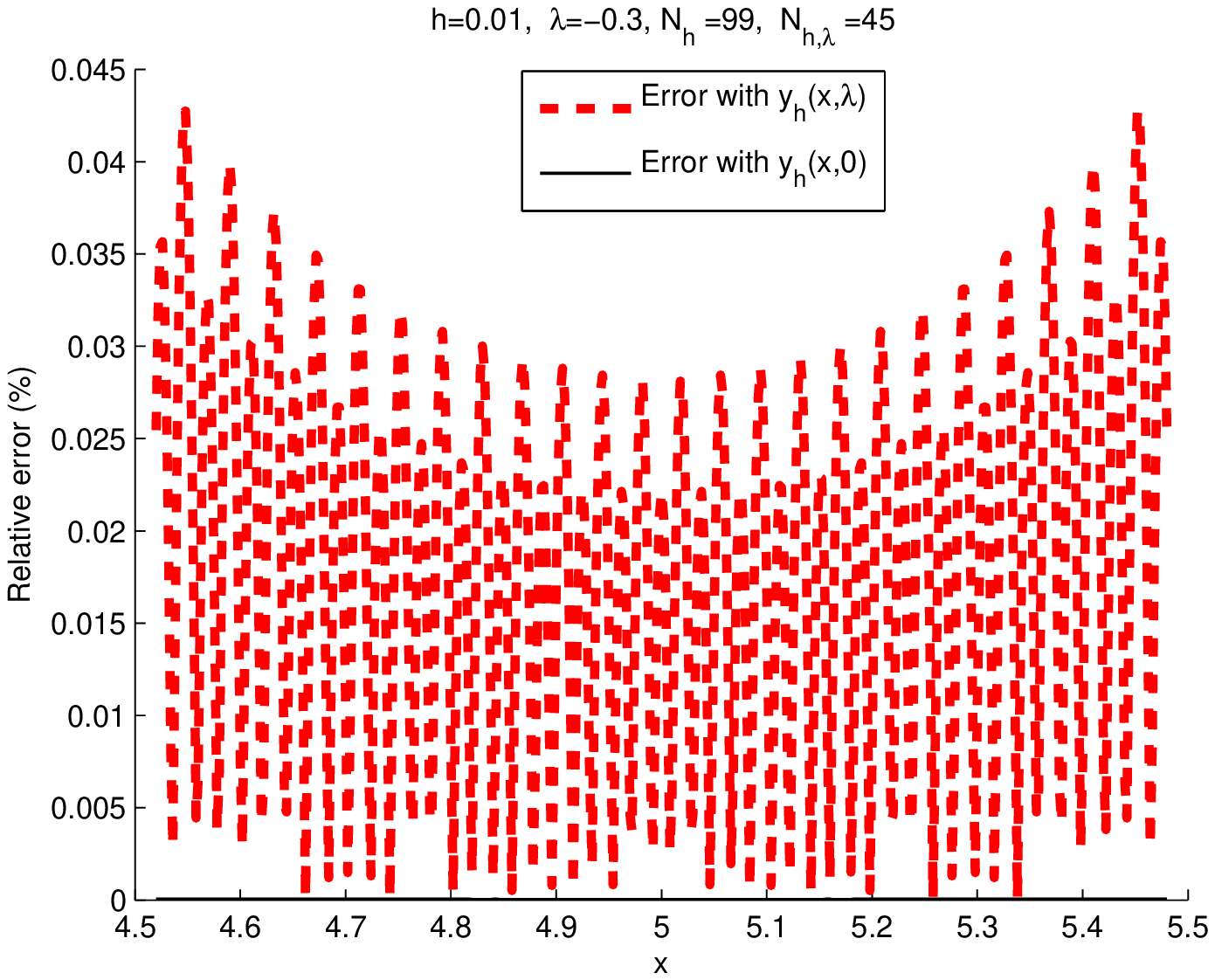}}
\subfigure[$h=0.001$]{\includegraphics[width=7cm]{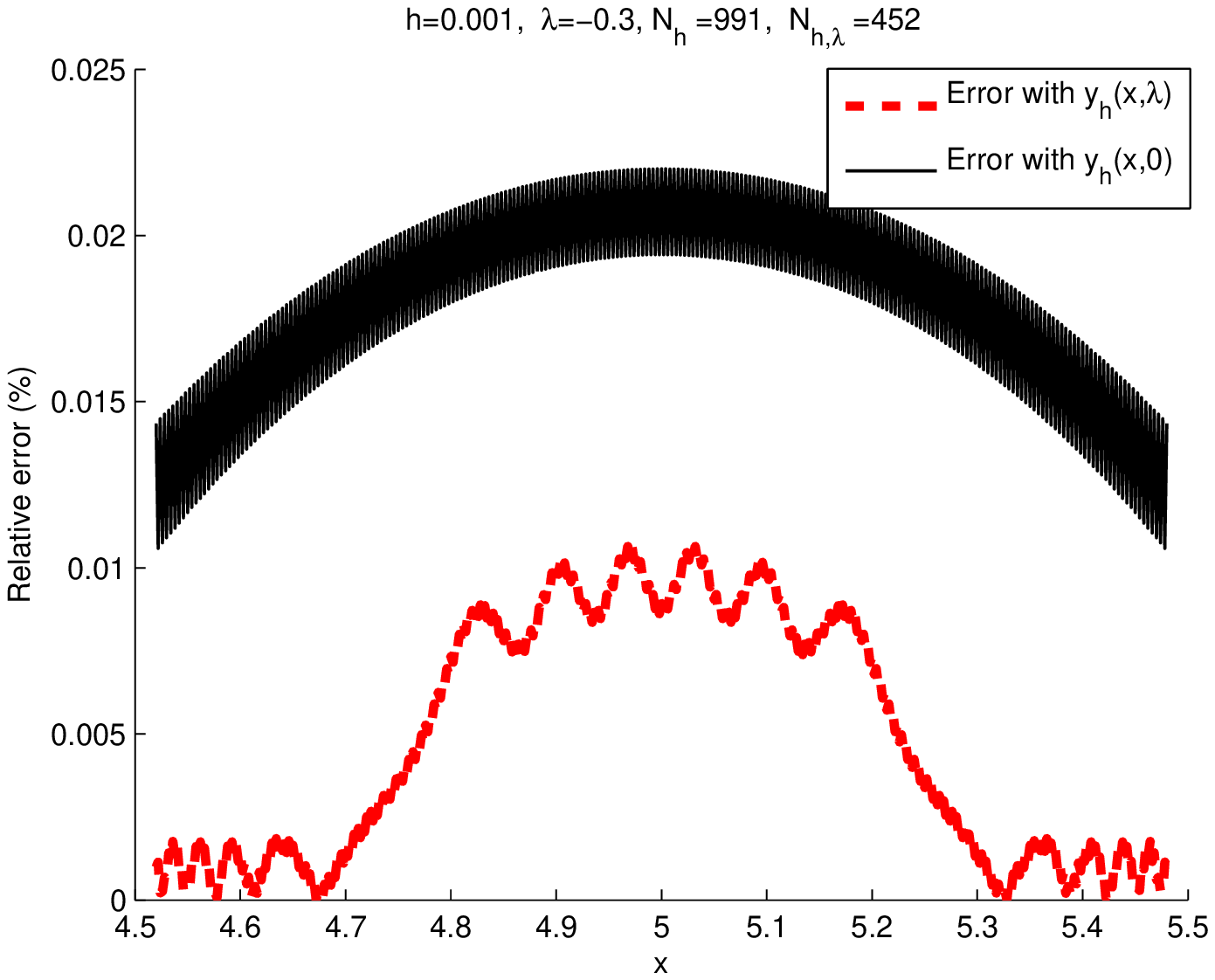}}
\caption{Relative error for $\lambda=0$ and $\lambda=-0.3$ for the $\sech^2$ example}
\label{erreur-sech-lamda03}
\end{center}
\end{figure}

\begin{figure}[htbp]
\begin{center}
\subfigure[$h=0.01$]{\includegraphics[width=7cm]{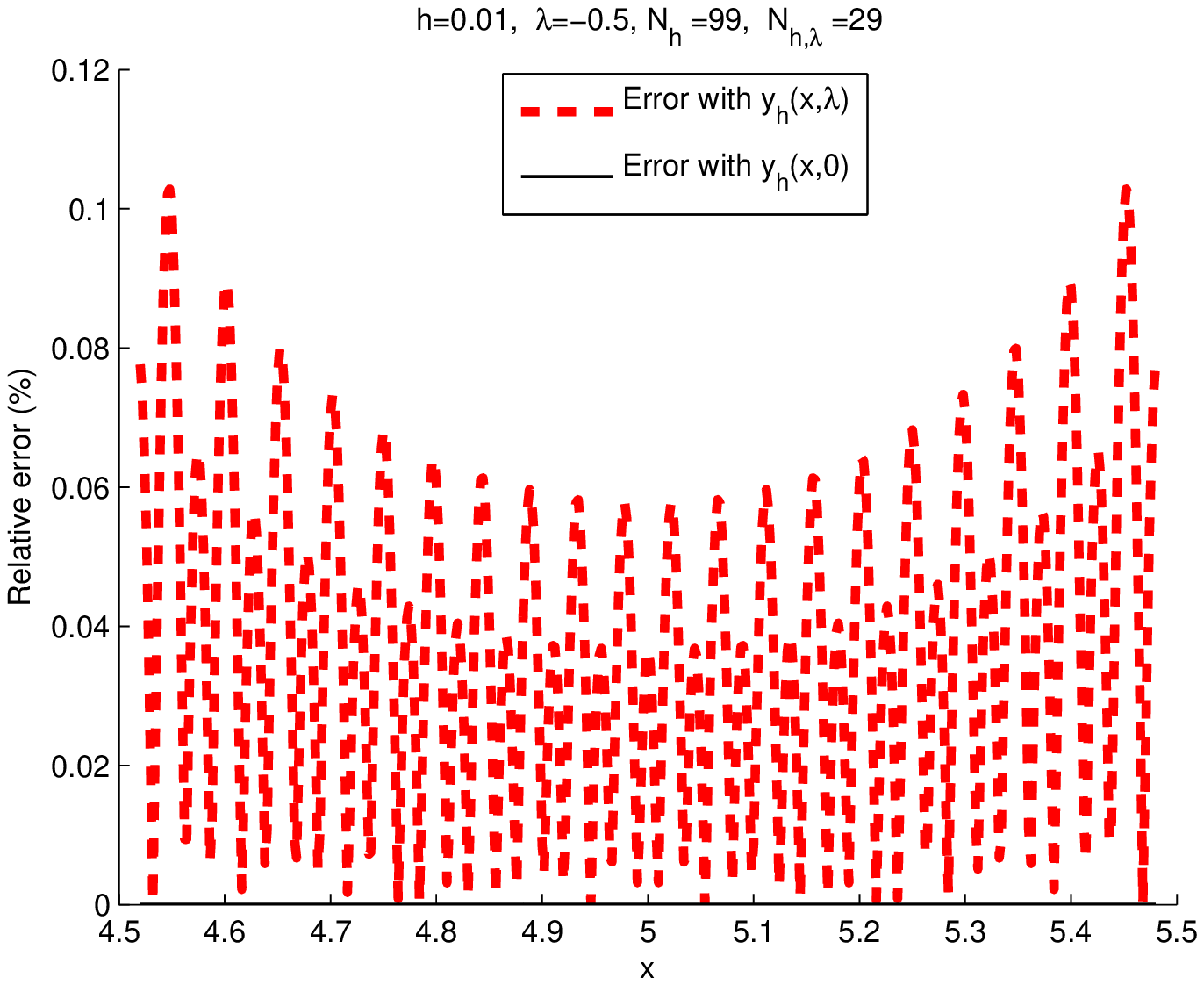}}
\subfigure[$h=0.001$]{\includegraphics[width=7cm]{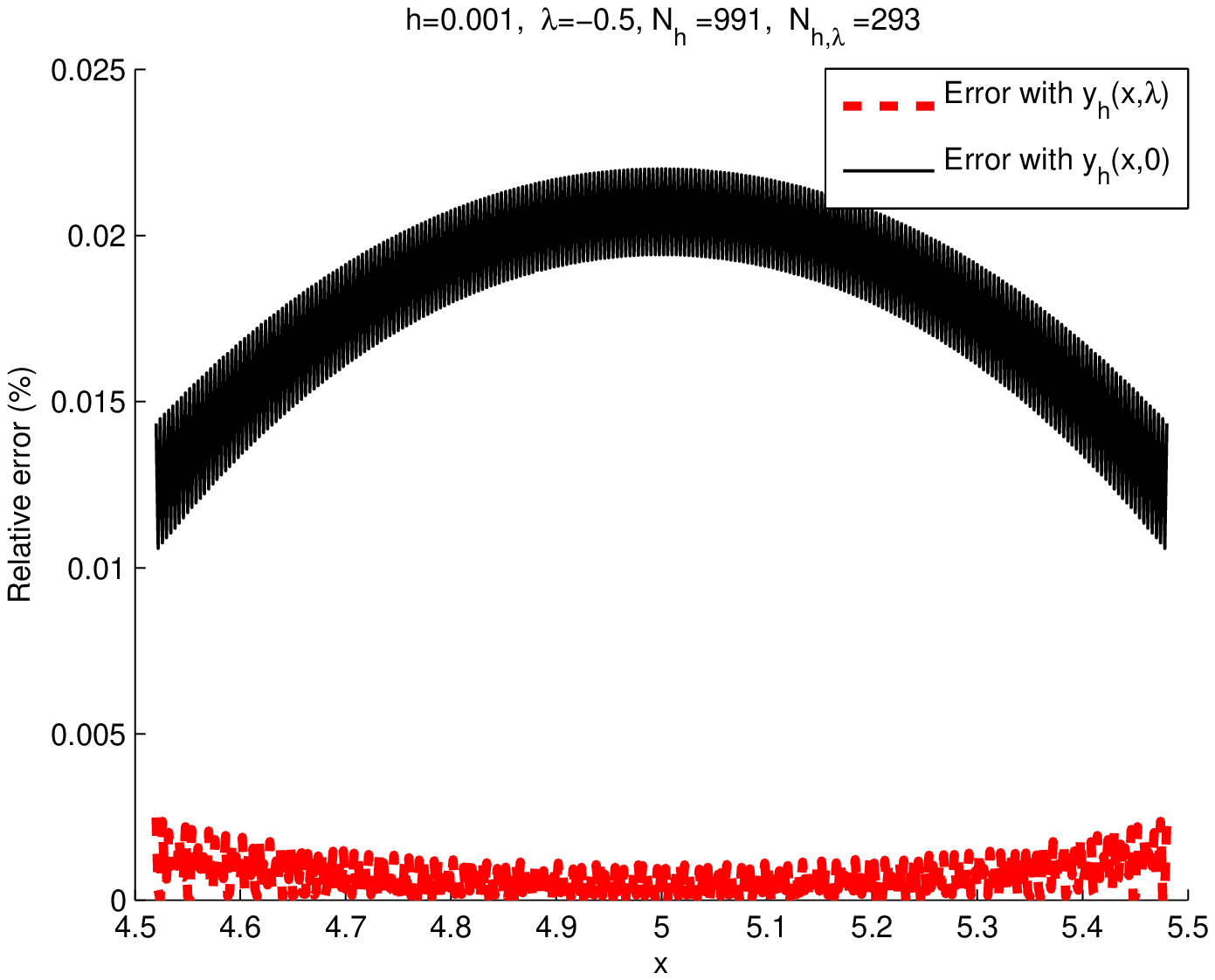}}
\caption{Relative error for $\lambda=0$ and $\lambda=-0.5$ for the $\sech^2$ example}
\label{erreur-sech-lamda05}
\end{center}
\end{figure}

\begin{figure}[htbp]
\begin{center}
\subfigure[]{\includegraphics[width=7cm]{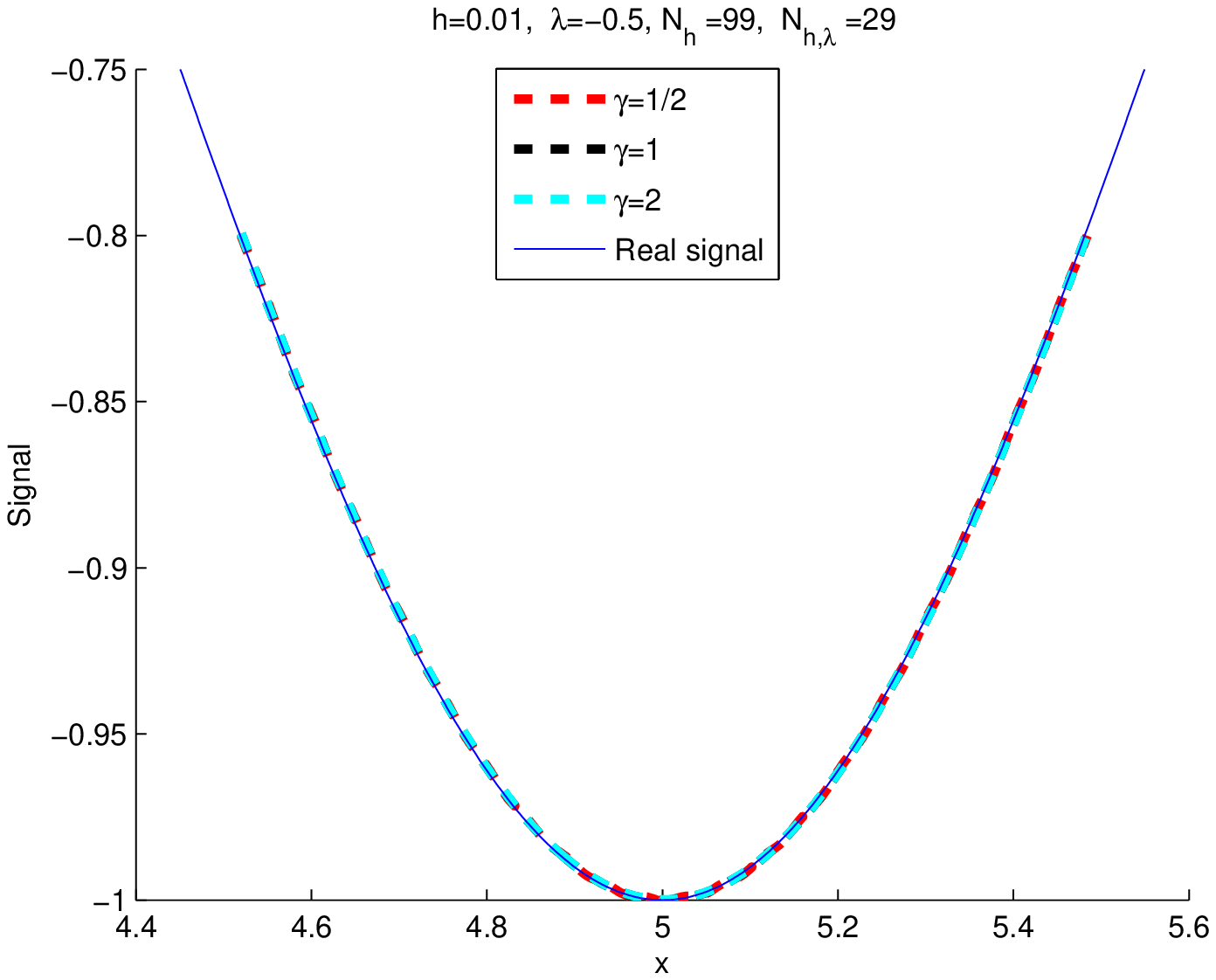}}
\subfigure[]{\includegraphics[width=7cm]{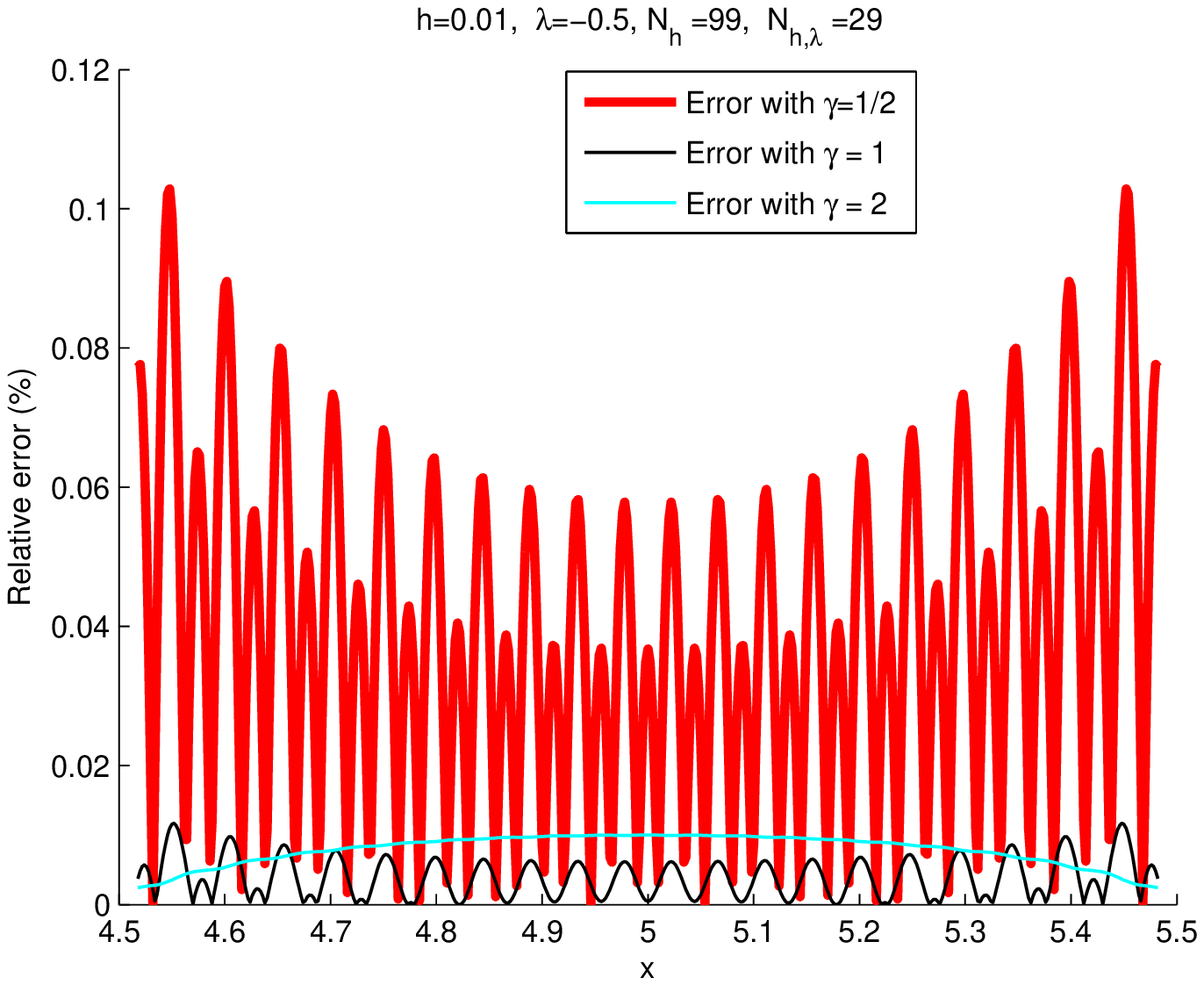}}
\caption{Reconstruction of a part of the $\sech^2$ signal (a) and relative error (b) for different values of $\gamma$ and $h=0.01$}
\label{h001-gamma05-1-2}
\end{center}
\end{figure}

\subsection{The arterial blood pressure (ABP) signal}
\begin{figure}[htbp]
\begin{center}
\subfigure[Example 1]{\includegraphics[width=7cm]{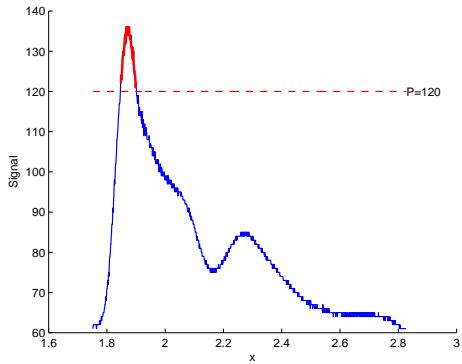}}
\subfigure[Example 2]{\includegraphics[width=7cm]{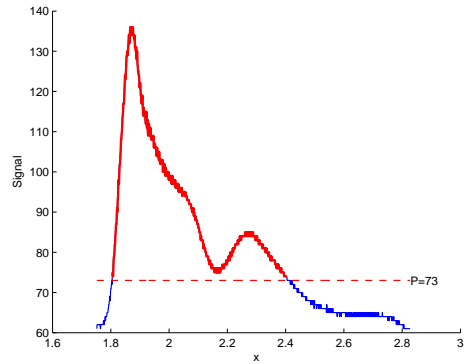}}
\caption{An example of an arterial blood pressure beat measured at the finger. The part of the signal that we want to reconstruct is in red.}
\label{pression-exemple}
\end{center}
\end{figure}

\begin{figure}[htbp]
\begin{center}
\subfigure[]{\includegraphics[width=7cm]{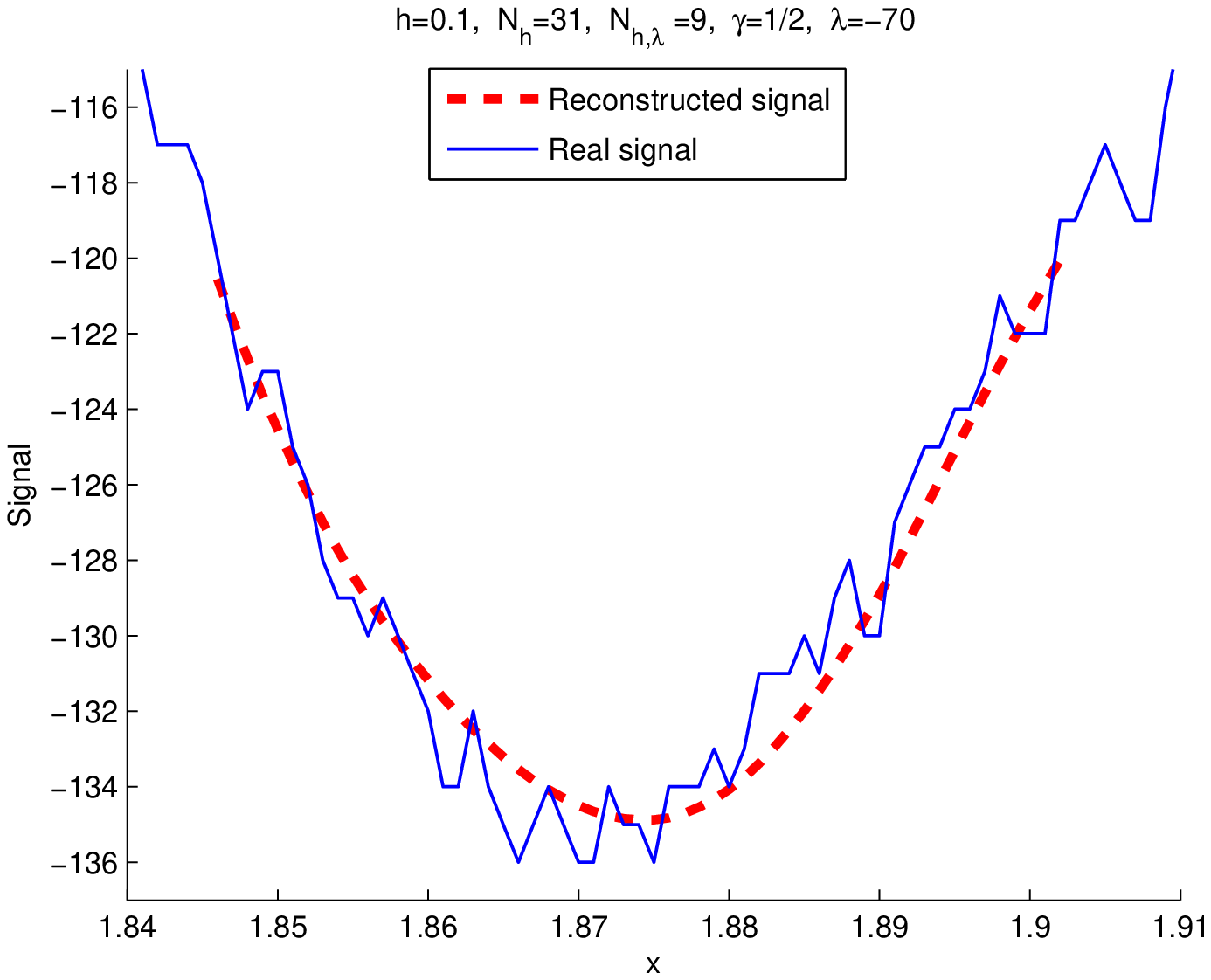}}
\subfigure[]{\includegraphics[width=7cm]{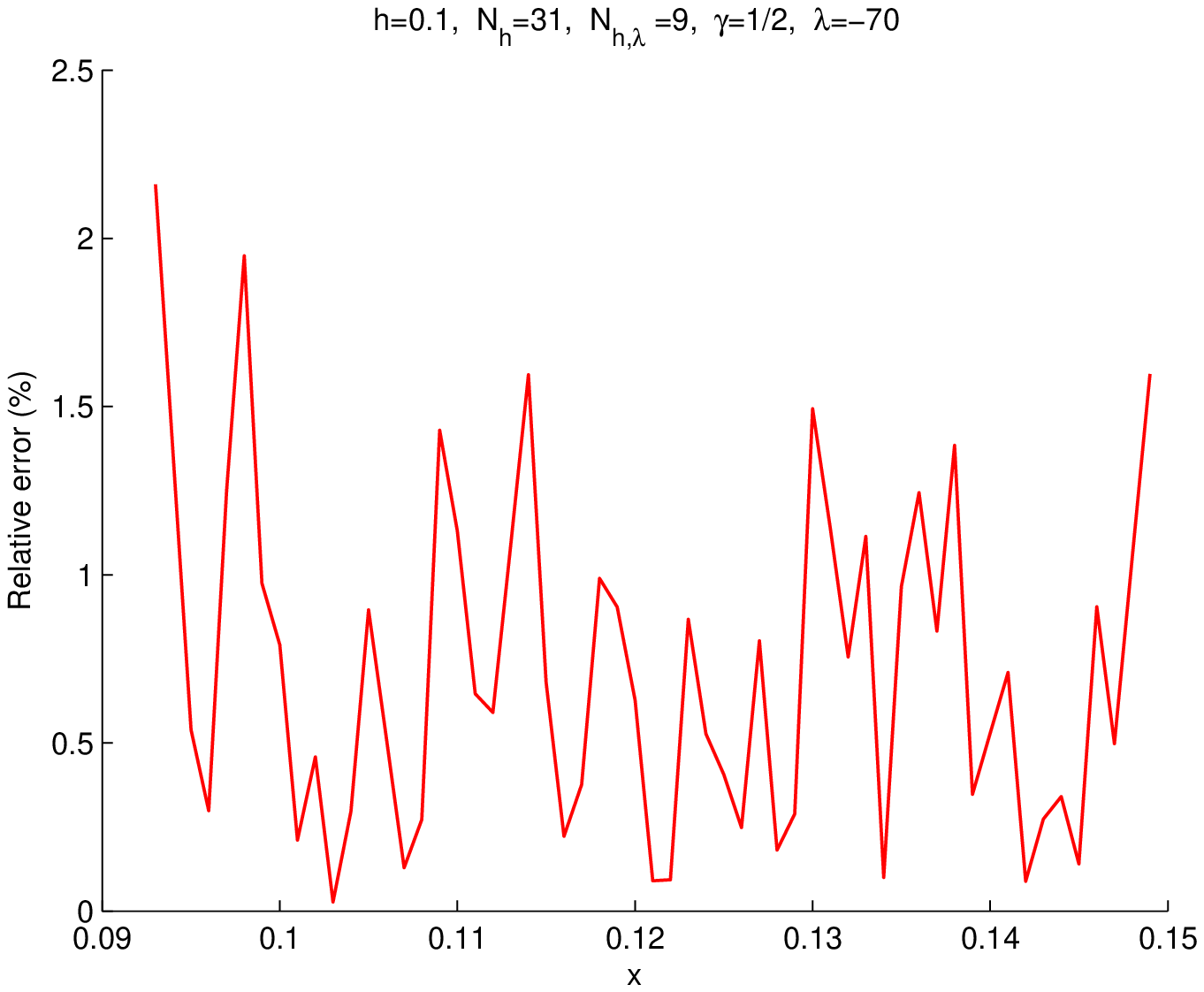}}
\caption{Reconstruction of a part of the ABP signal (a)  and relative error (b) for $\lambda=-70$, $h=0.1$ and $\gamma=\frac{1}{2}$}
\label{pression-h01-lamda70-gamma05}
\end{center}
\end{figure}

\begin{figure}[htbp]
\begin{center}
\subfigure[]{\includegraphics[width=7cm]{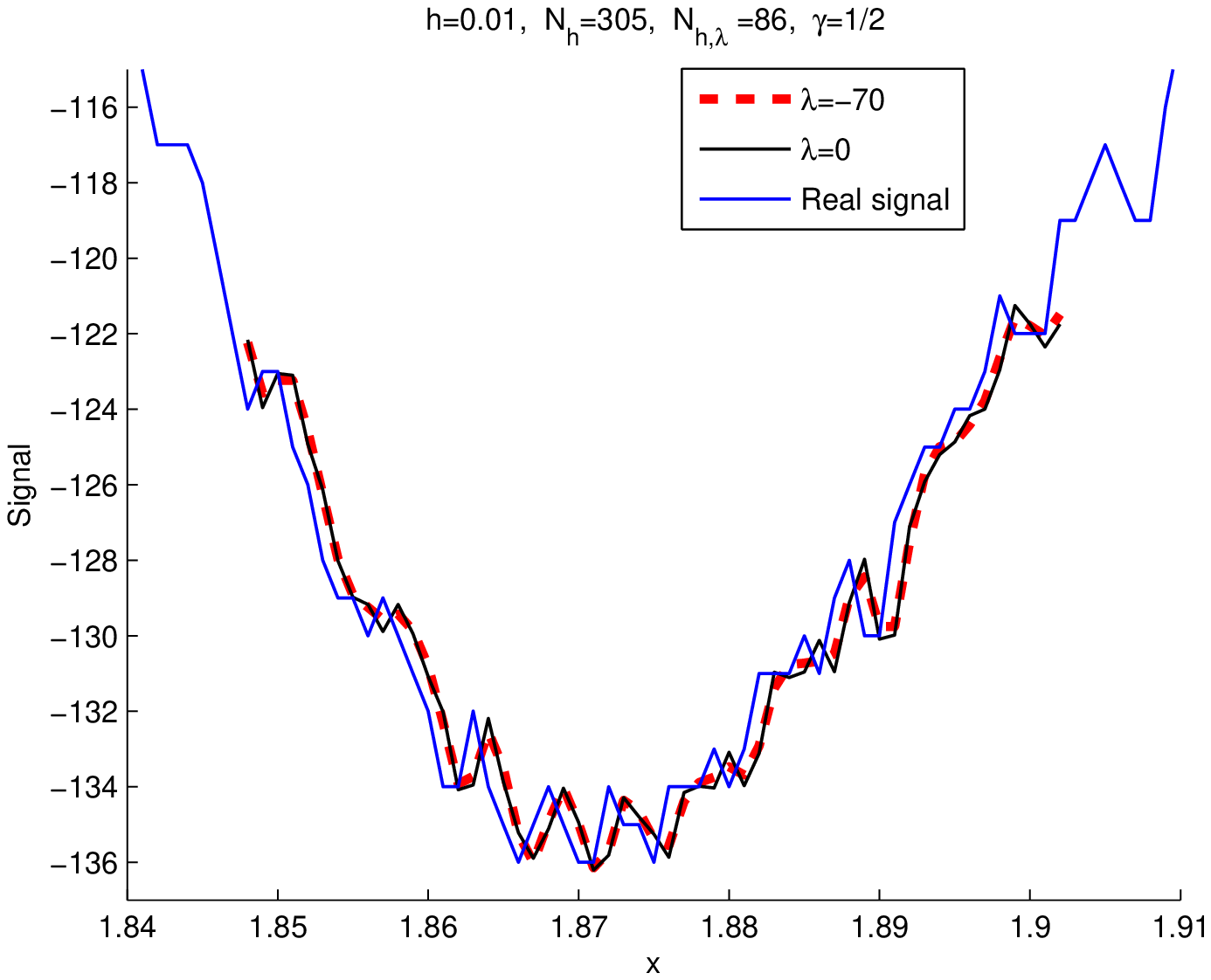}}
\subfigure[]{\includegraphics[width=7cm]{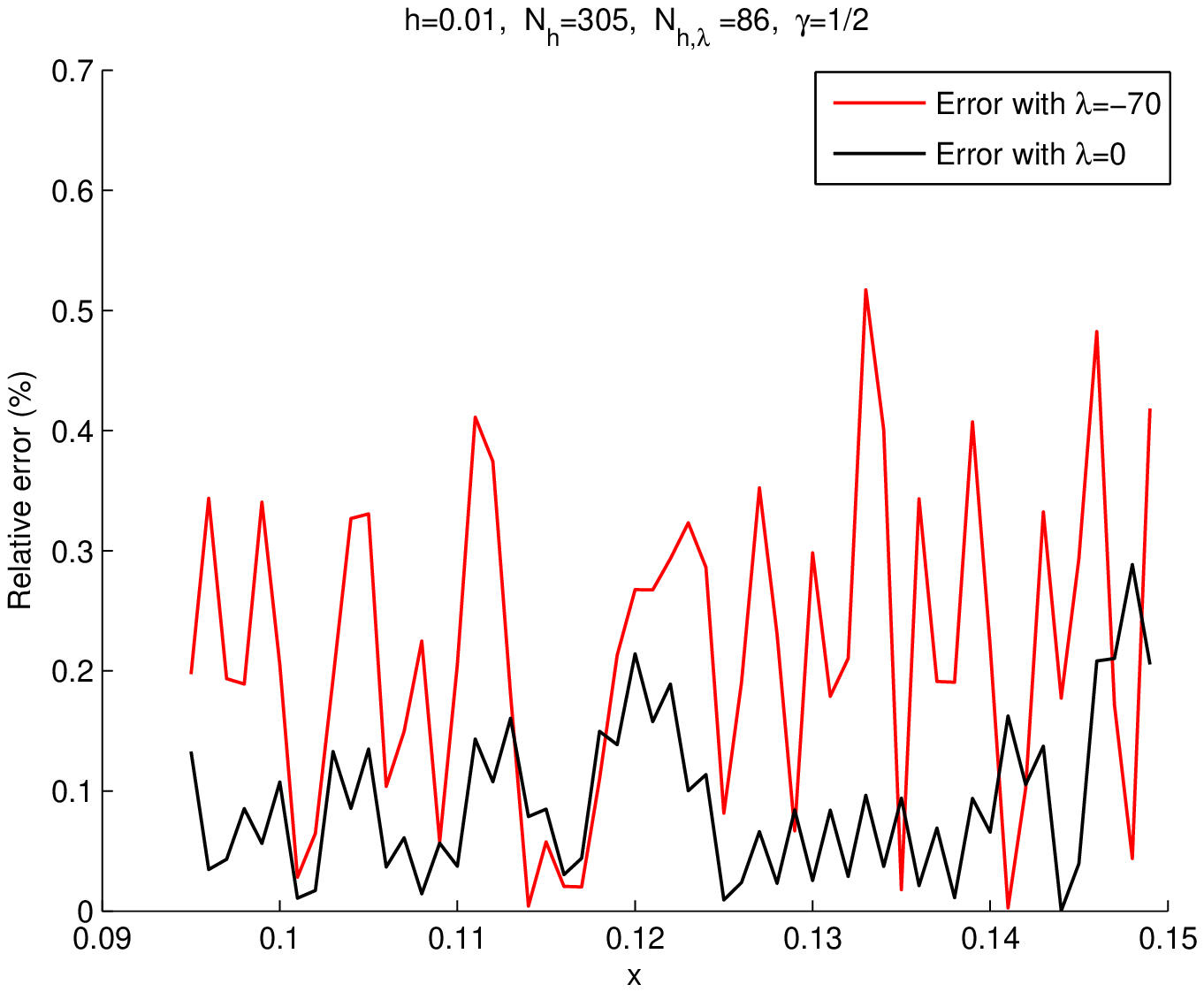}}
\caption{Reconstruction of a part of the ABP signal (a) and relative error (b) for $\lambda=0$ and $-70$, $h=0.01$ and $\gamma=\frac{1}{2}$}
\label{pression-h001-lamda0-70-gamma05}
\end{center}
\end{figure}

\begin{figure}[htbp]
\begin{center}
\subfigure[]{\includegraphics[width=7cm]{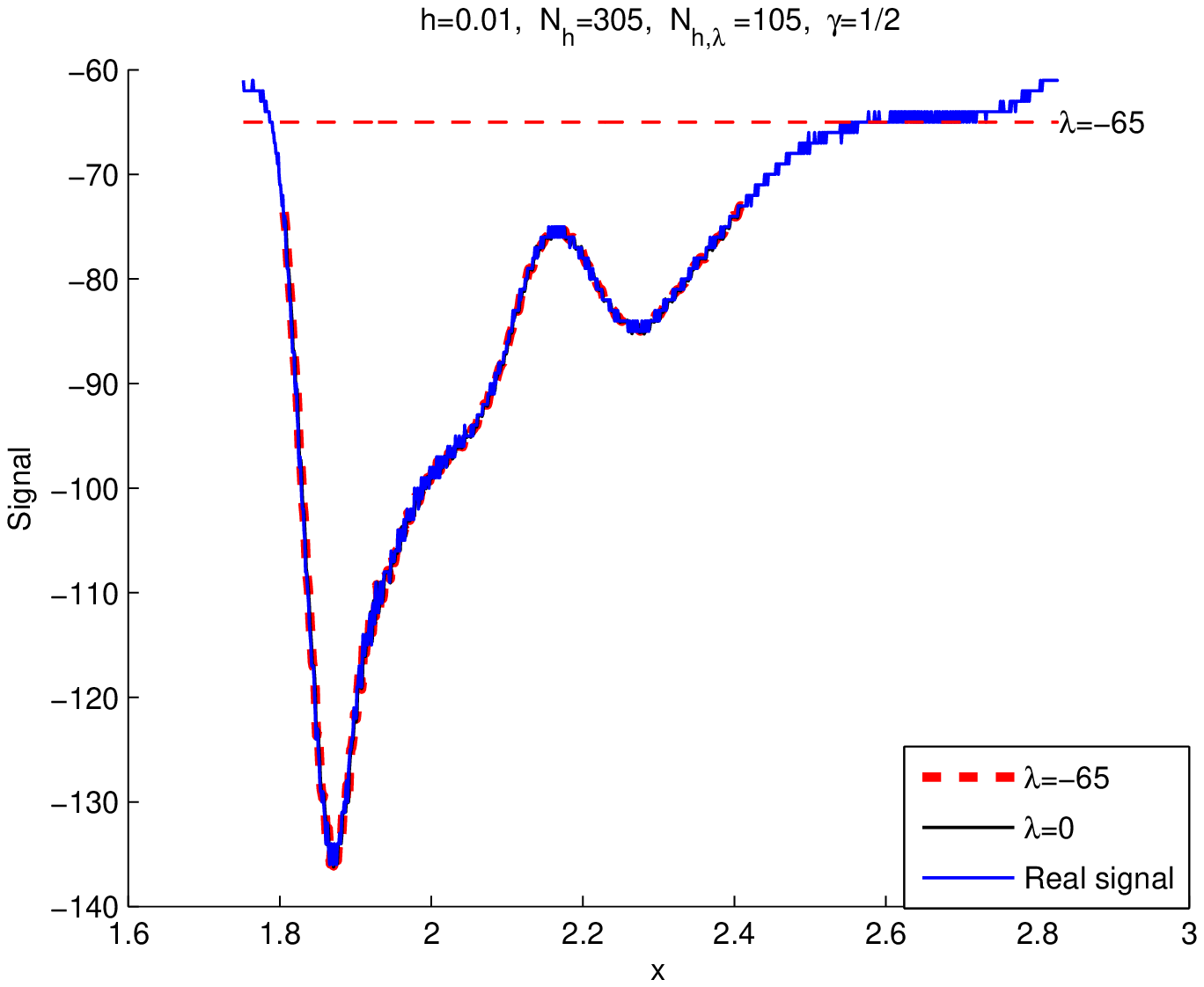}}
\subfigure[]{\includegraphics[width=7cm]{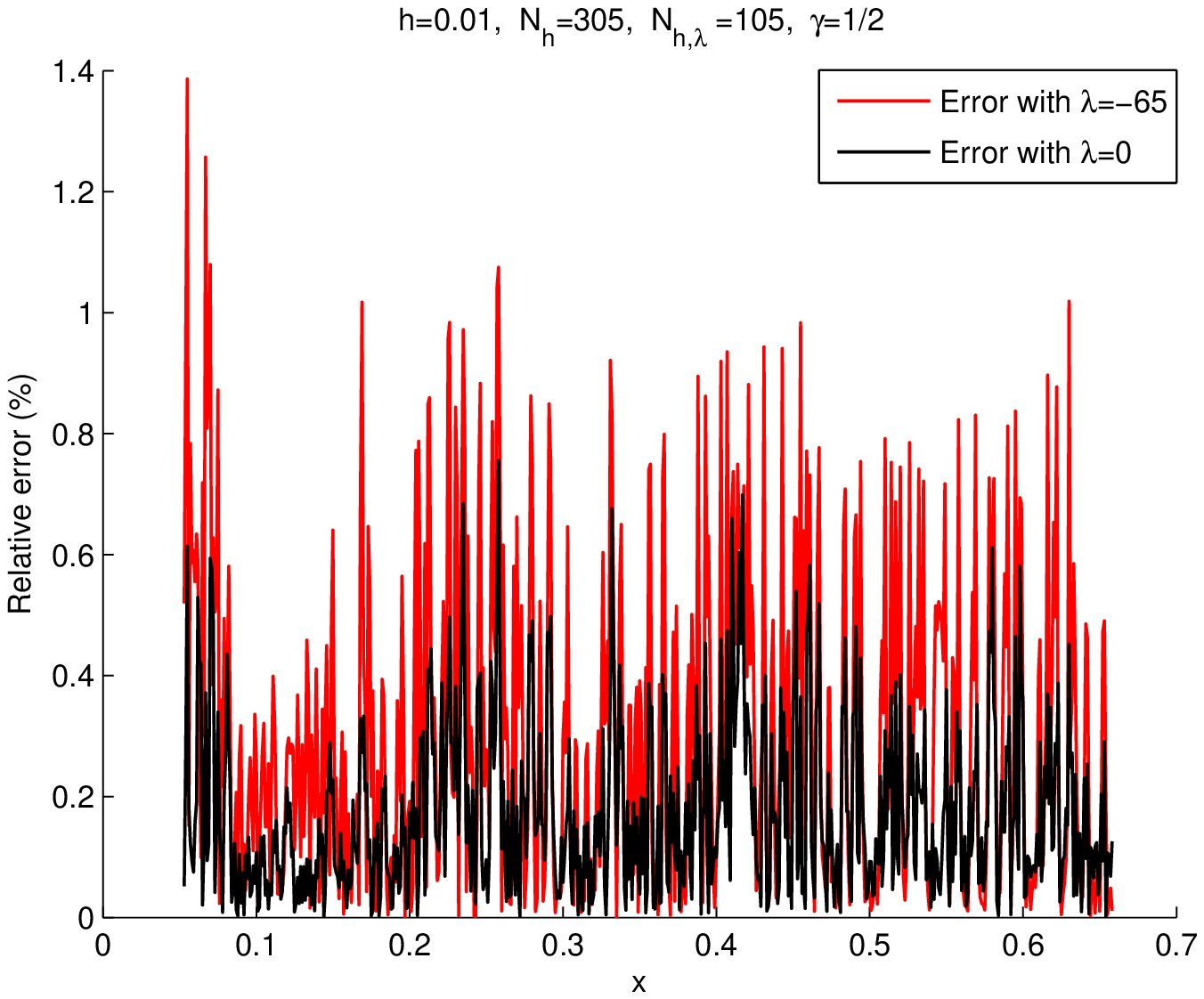}}
\caption{Reconstruction of a part of the ABP signal (a) and relative error (b) for $\lambda=0$ and $-65$, $h=0.01$ and $\gamma=\frac{1}{2}$}
\label{pression-h001-lamda0-65-gamma05}
\end{center}
\end{figure}

We analyze in this section the results obtained when the signal $y(x)$ is an ABP signal. ABP  plays an important role in the cardiovascular system and is used in clinical practice for monitoring purposes. However the interpretation of ABP signals is still restricted to the interpretation of the maximal and the minimal values called respectively the systolic and diastolic pressures.
No information on the instantaneous variability of the pressure is considered. Recent studies have proposed to exploit the ABP waveform in clinical practice by analyzing the signal with a semi-classical signal analysis approach (see for example \cite{LaCrSo:10}, \cite{LaMePaCoVa:10}). The latter proposes to reconstruct the signal with formula (\ref{formule2}), corresponding to the case $\lambda=0$ in order to extract some spectral quantities that provide an  interesting information on that signal. These quantities are the eigenvalues and some Riesz means of these eigenvalues.  These quantities enable for example the discrimination between different pathological and physiological situations \cite{LaMeCoSo:07} and also provide information on some cardiovascular parameters of great interest as for example the stroke volume \cite{LaMePaCoVa:10}. We are interested in this study in the reconstruction  of a small part of an ABP pressure beat illustrated in fig. \ref{pression-exemple} with formula \eqref{formule}. We consider different values of $h$, $\lambda$ and $\gamma$ as it is described in fig. \ref{pression-h01-lamda70-gamma05} - fig. \ref{pression-h001-lamda0-65-gamma05} where $-y$ is represented. In fig. \ref{pression-h01-lamda70-gamma05} - fig. \ref{pression-h001-lamda0-70-gamma05} a zoom on the reconstructed part of the signal is represented. Note that the signal in this case is not a regular function but only known from data measurements for a sequence $x_j = a + (j-1)\dfrac{b-a}{M-1}$ ($j=1,\dots,M$) for some integer $M$.   We recall that in our application the $x$ variable represents the time. The time between two consecutive measurements is $10^{-3}$ seconds.

\begin{figure}[htbp]
\begin{center}
\subfigure[]{\includegraphics[width=7cm]{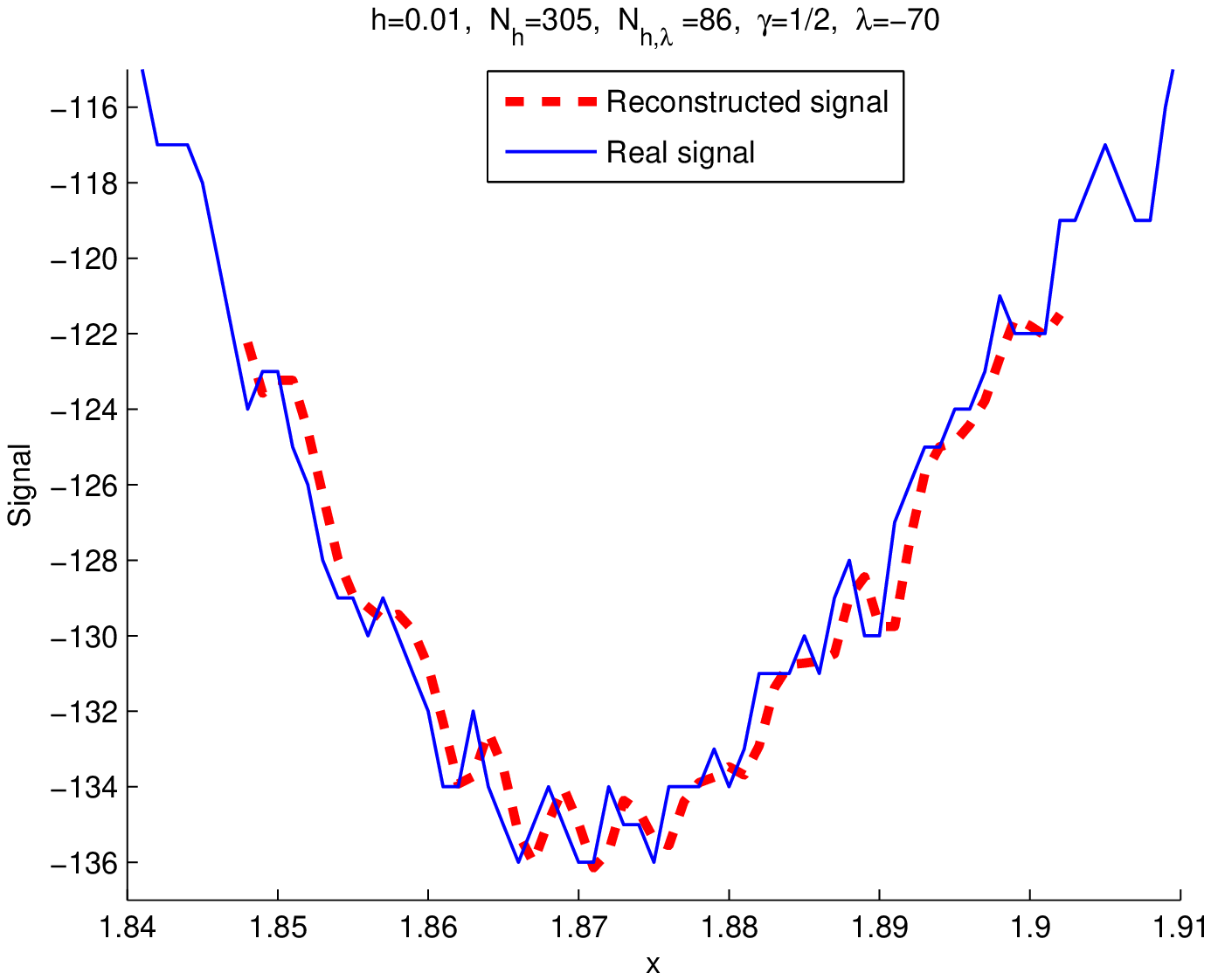}}
\subfigure[]{\includegraphics[width=7cm]{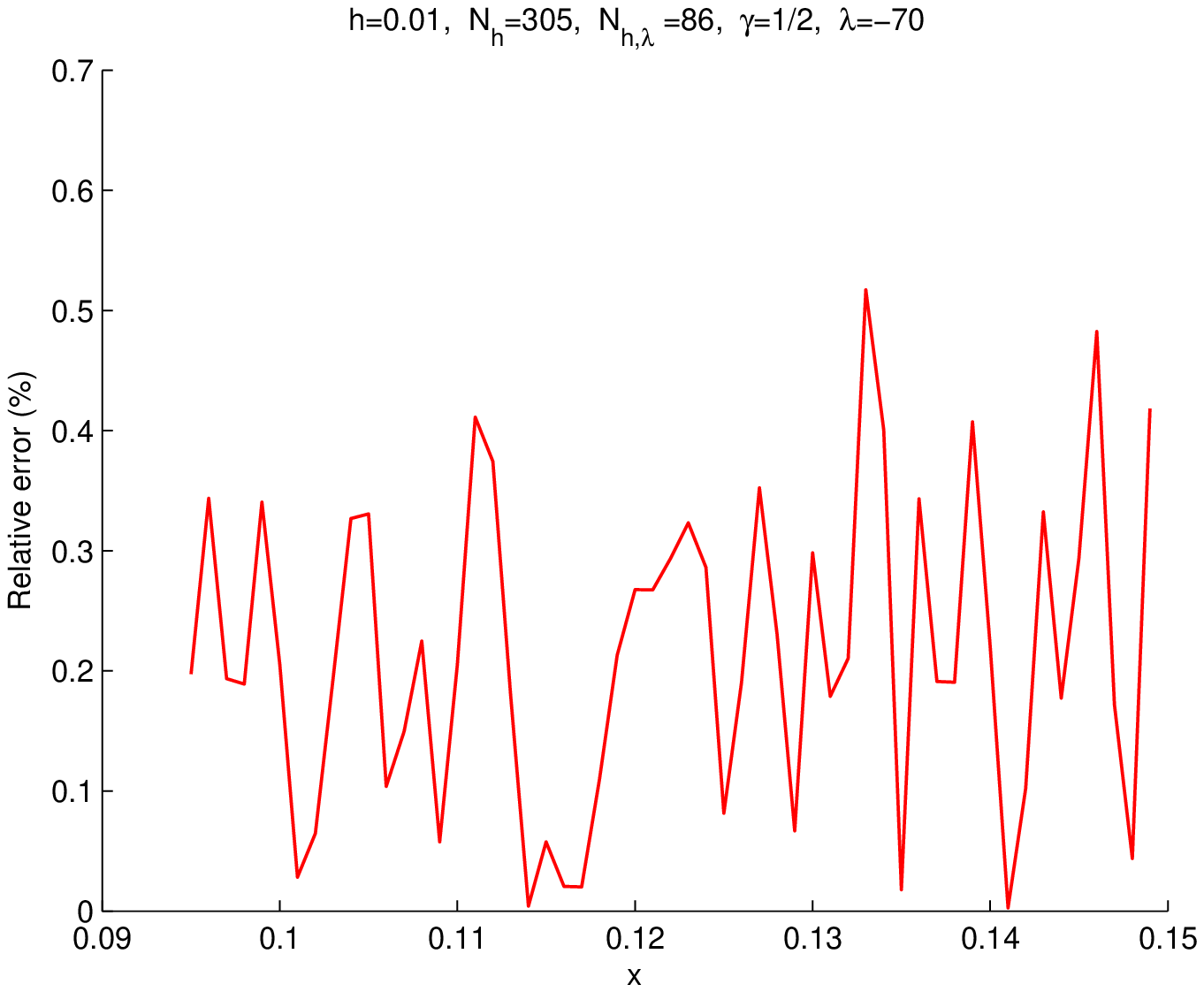}}
\caption{Reconstruction of a part of the ABP signal (a) and relative error (b) for $\lambda=-70$, $h=0.01$ and $\gamma=\frac{1}{2}$}
\label{pression-h001-lamda70-gamma05}
\end{center}
\end{figure}

\begin{figure}[htbp]
\begin{center}
\subfigure[]{\includegraphics[width=7cm]{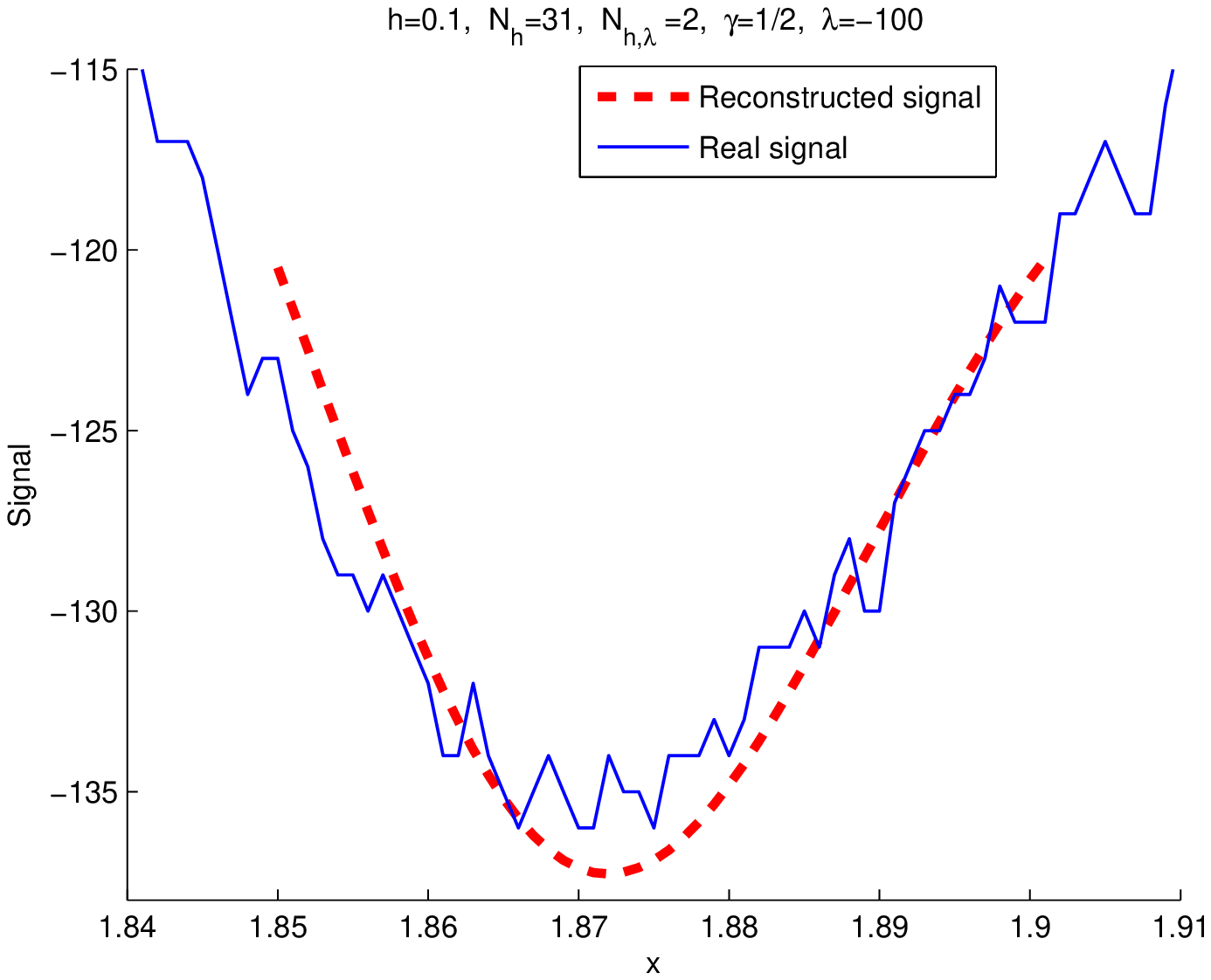}}
\subfigure[]{\includegraphics[width=7cm]{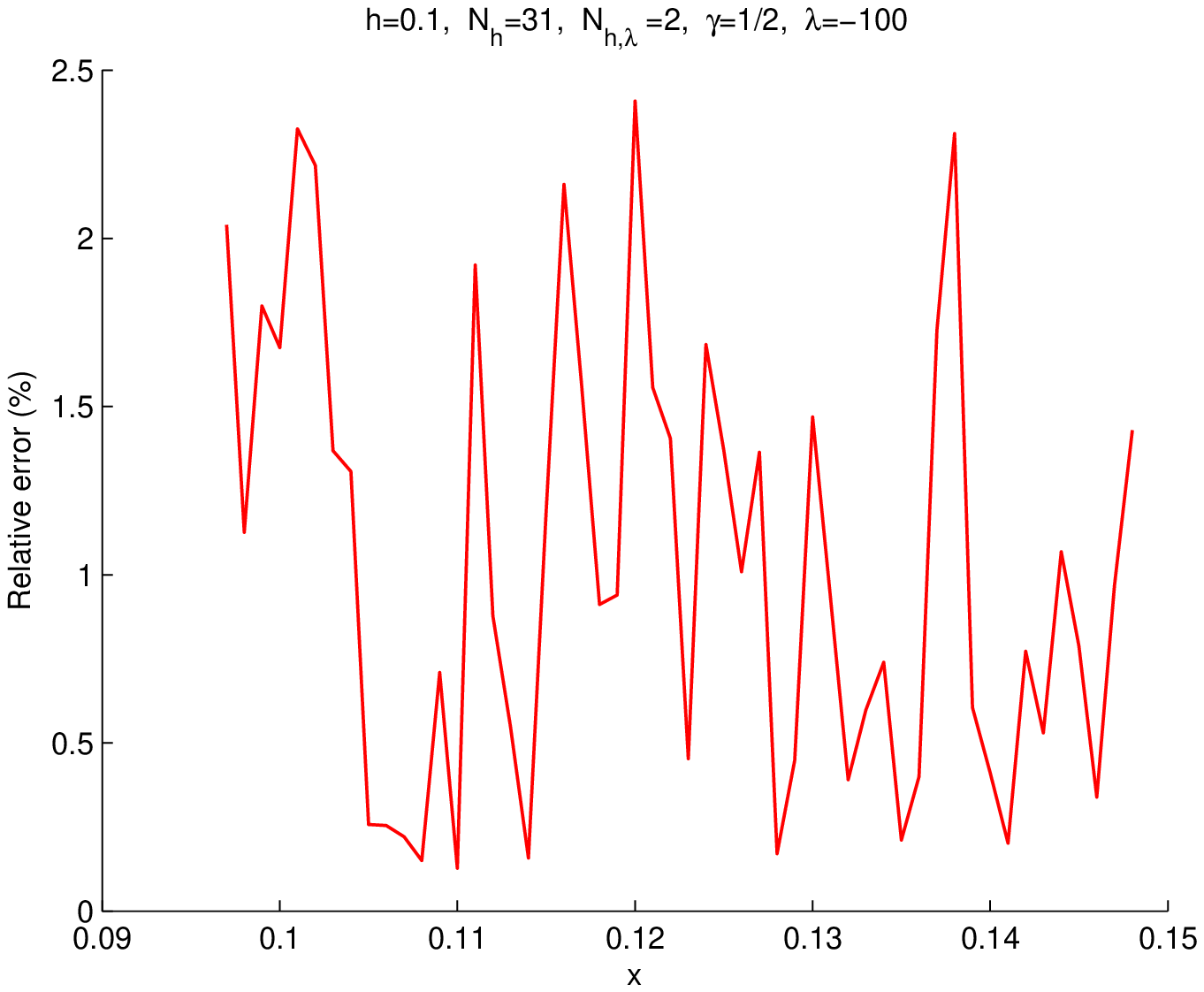}}
\caption{Reconstruction of a part of the ABP signal (a) and relative error (b) for $\lambda=-100$, $h=0.1$ and $\gamma=\frac{1}{2}$}
\label{pression-h01-lamda100-gamma05}
\end{center}
\end{figure}

\begin{figure}[htbp]
\begin{center}
\subfigure[]{\includegraphics[width=7cm]{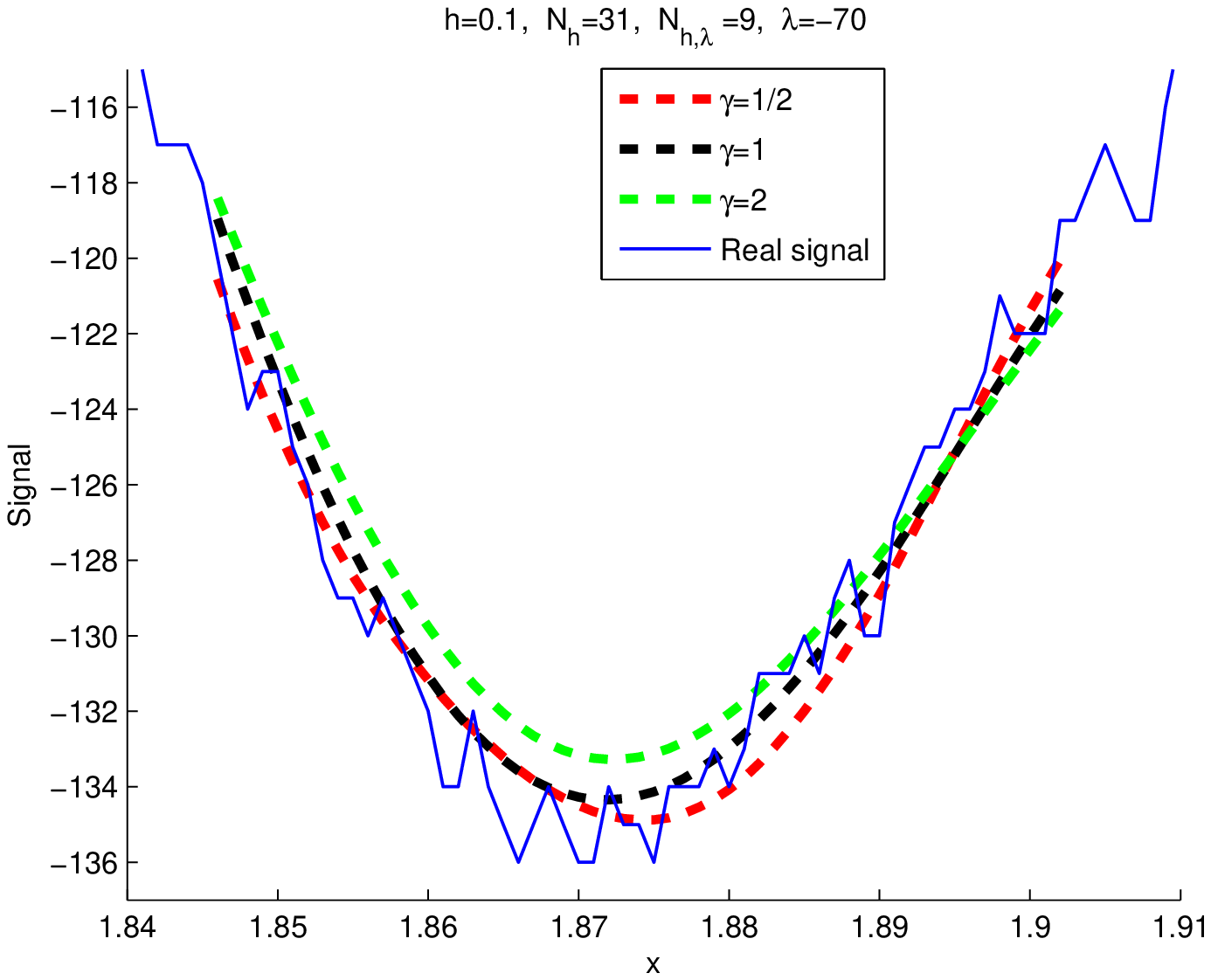}}
\subfigure[]{\includegraphics[width=7cm]{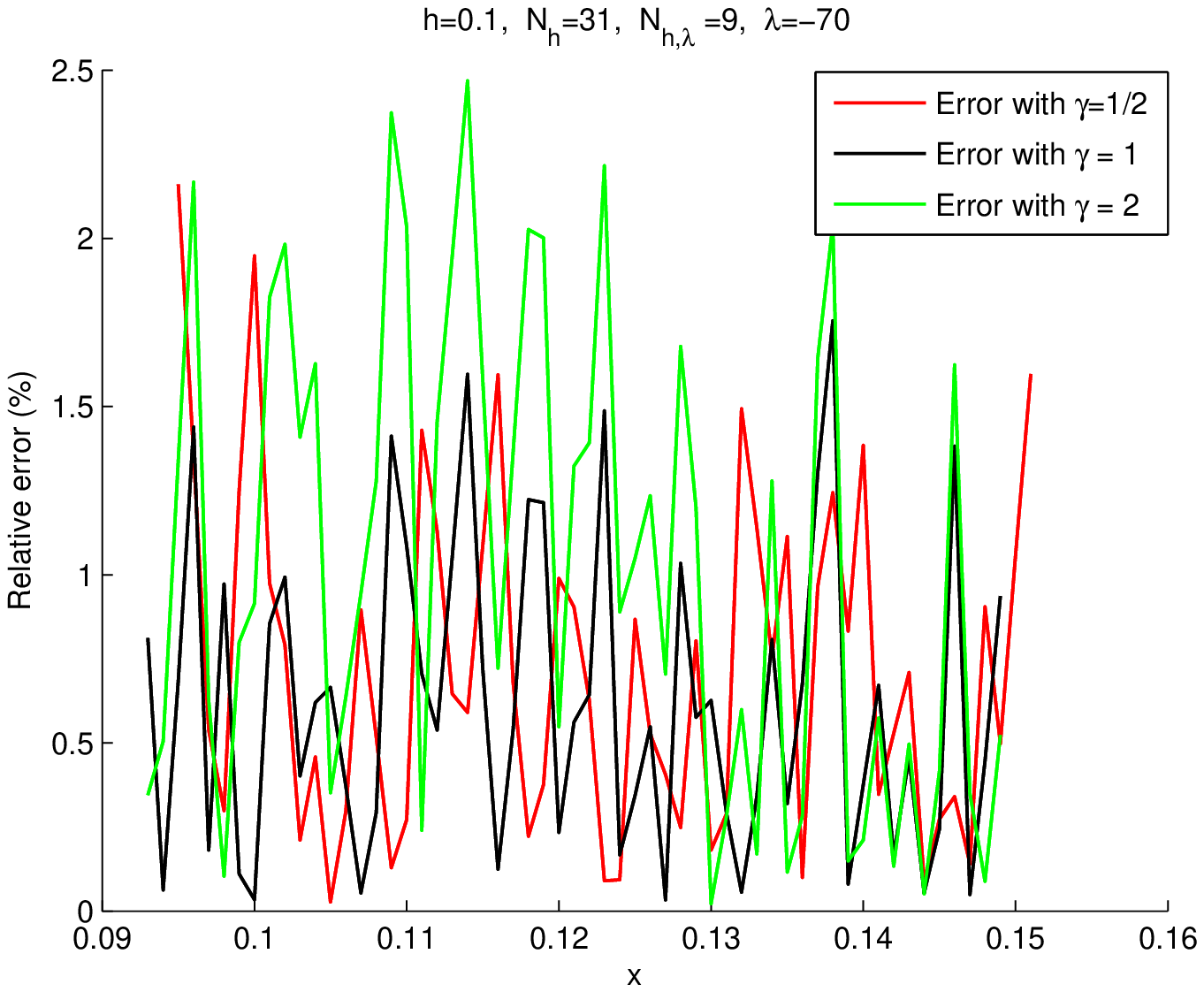}}
\caption{Reconstruction of a part of the ABP signal (a) and relative error (b) for $\lambda=-70$, $h=0.1$ and $\gamma=\frac{1}{2}, 1, 2$}
\label{pression-h01-lamda70-gamma05-1-2}
\end{center}
\end{figure}

\begin{figure}[htbp]
\begin{center}
\subfigure[]{\includegraphics[width=7cm]{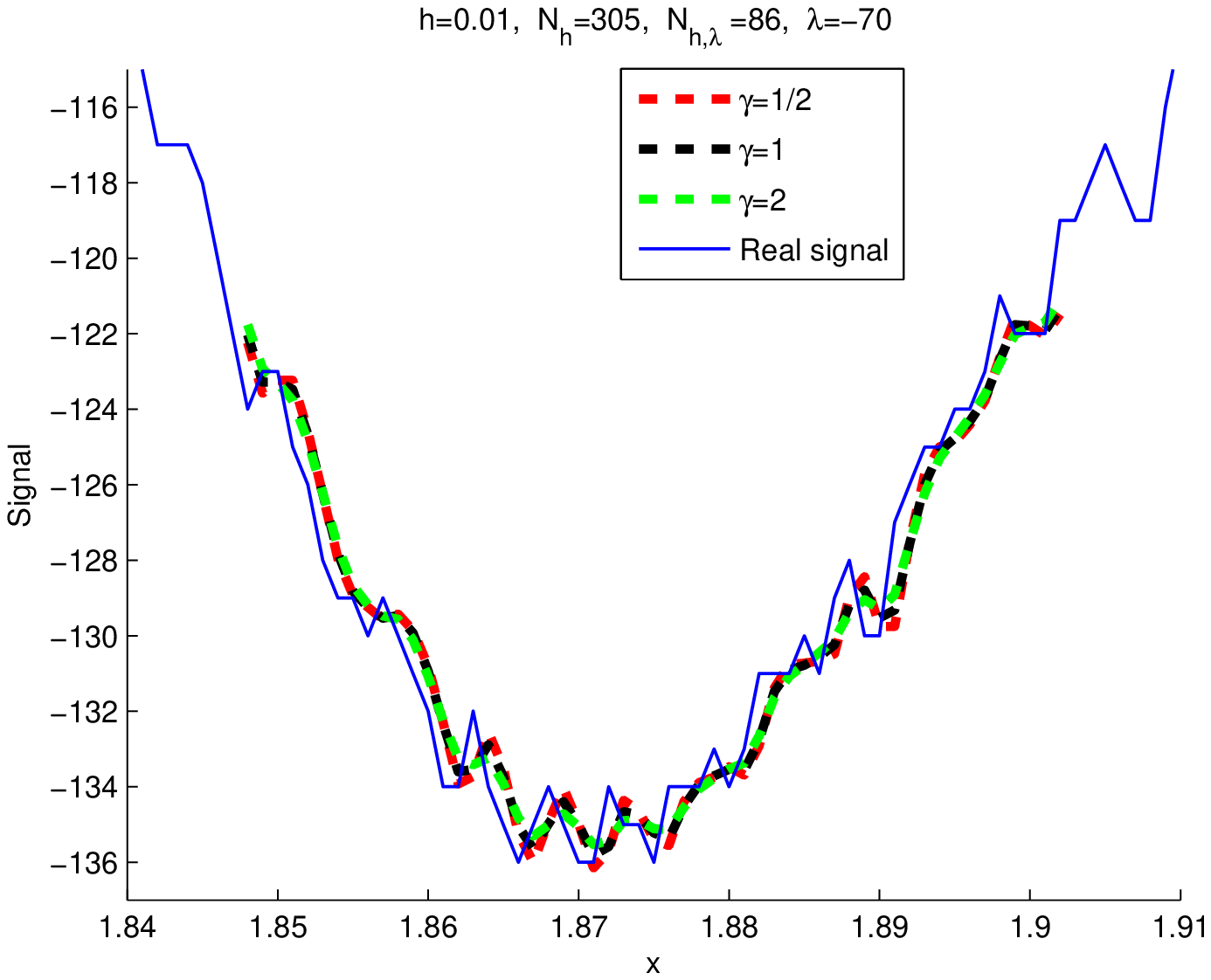}}
\subfigure[]{\includegraphics[width=7cm]{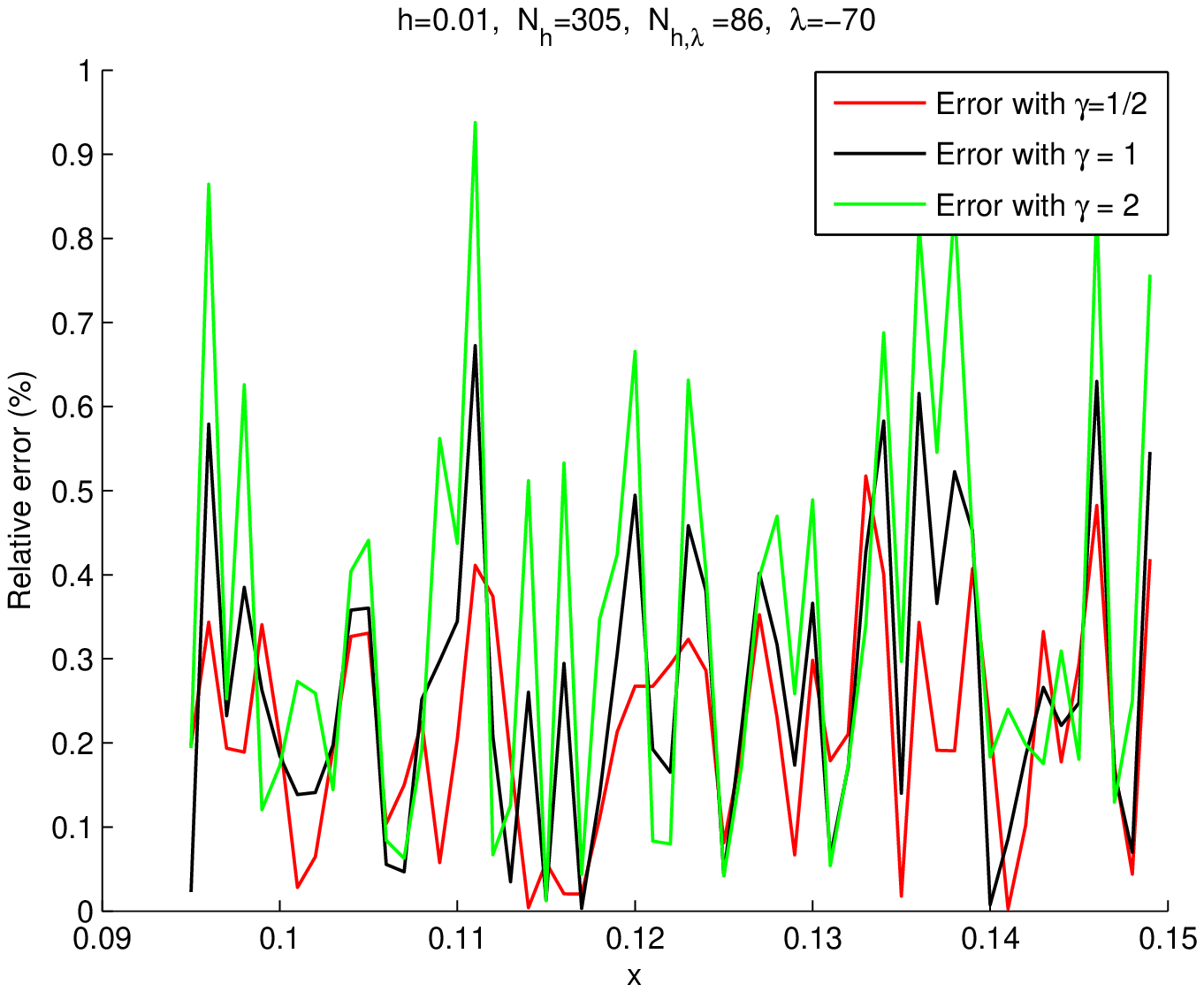}}
\caption{Reconstruction of a part of the ABP signal (a) and relative error (b) for $\lambda=-70$, $h=0.01$ and $\gamma=\frac{1}{2}, 1, 2$}
\label{pression-h001-lamda70-gamma05-1-2}
\end{center}
\end{figure}

\begin{figure}[htbp]
\begin{center}
\subfigure[]{\includegraphics[width=7cm]{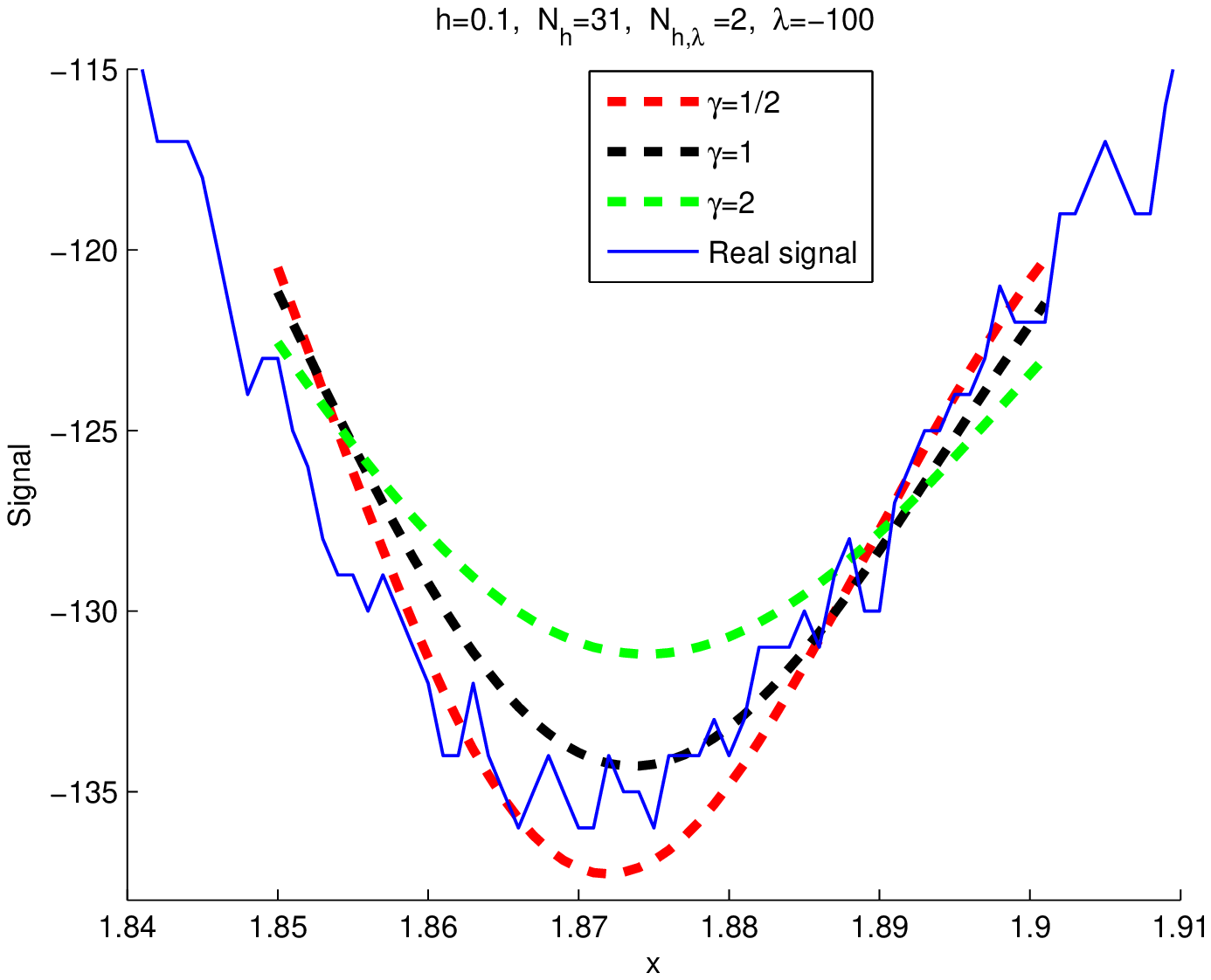}}
\subfigure[]{\includegraphics[width=7cm]{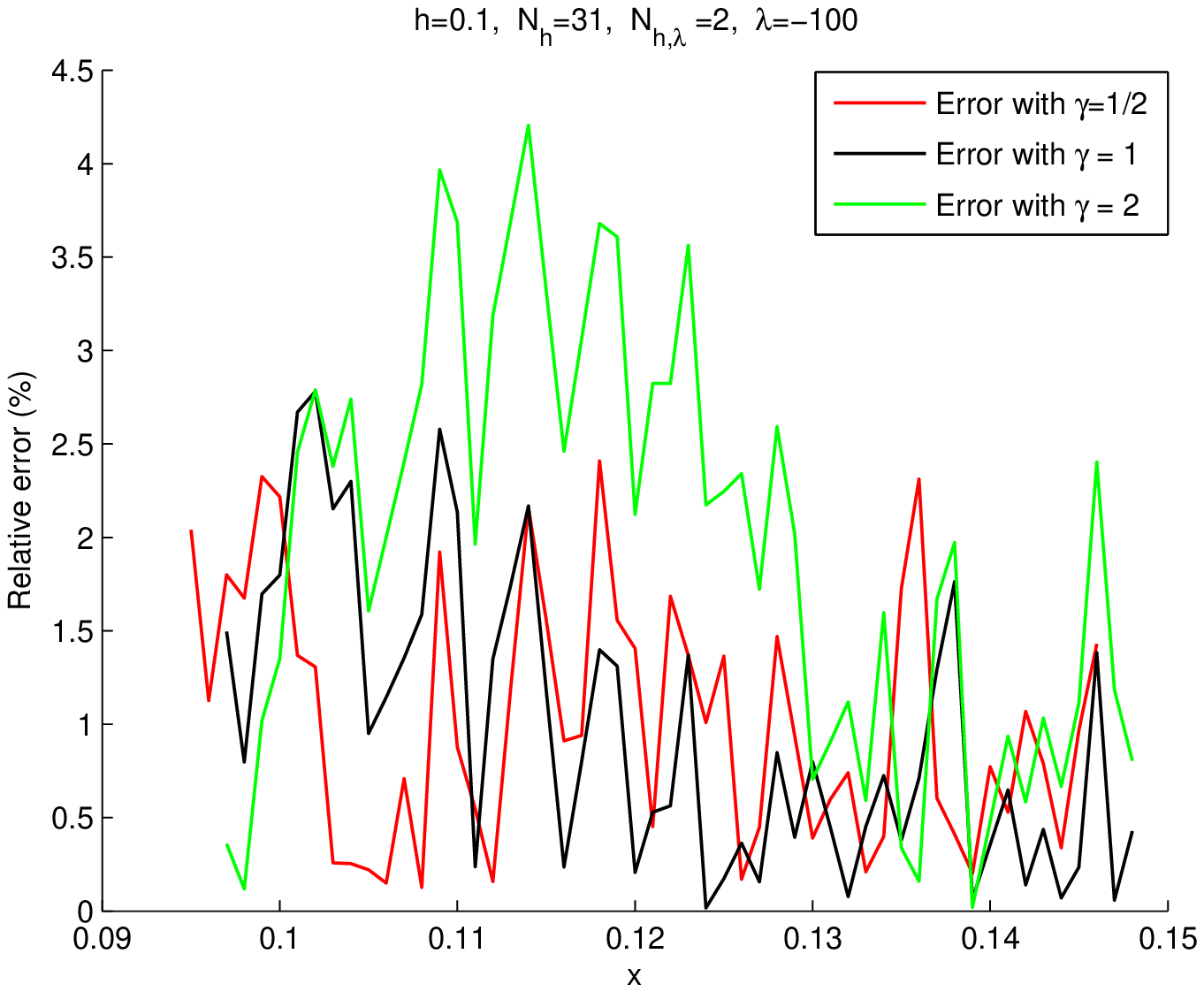}}
\caption{Reconstruction of a part of the ABP  signal (a) and relative error (b) for $\lambda=-100$, $h=0.1$ and $\gamma=\frac{1}{2}, 1, 2$}
\label{pression-h01-lamda100-gamma05-1-2}
\end{center}
\end{figure}

\begin{figure}[htbp]
\begin{center}
\subfigure[]{\includegraphics[width=7cm]{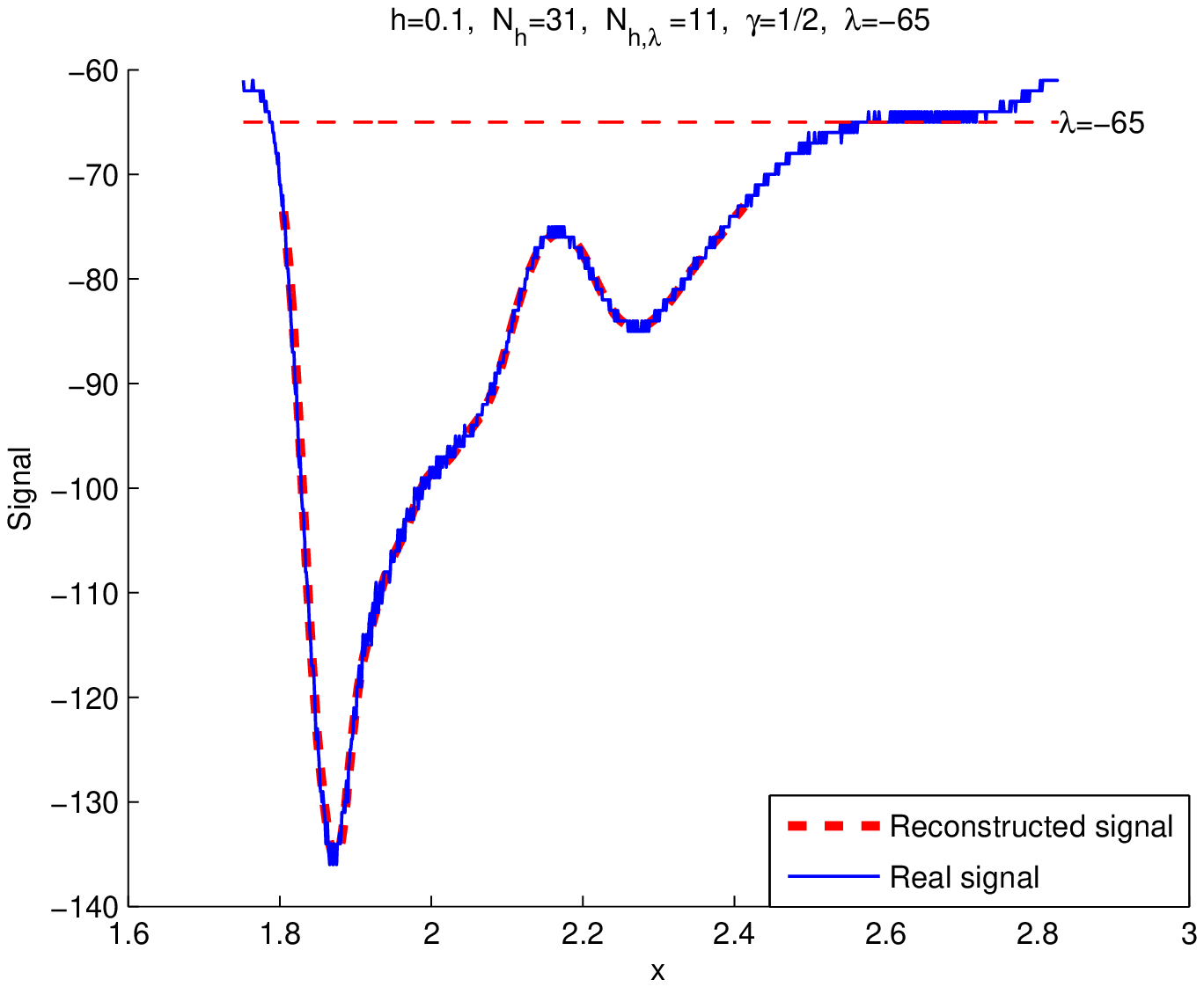}}
\subfigure[]{\includegraphics[width=7cm]{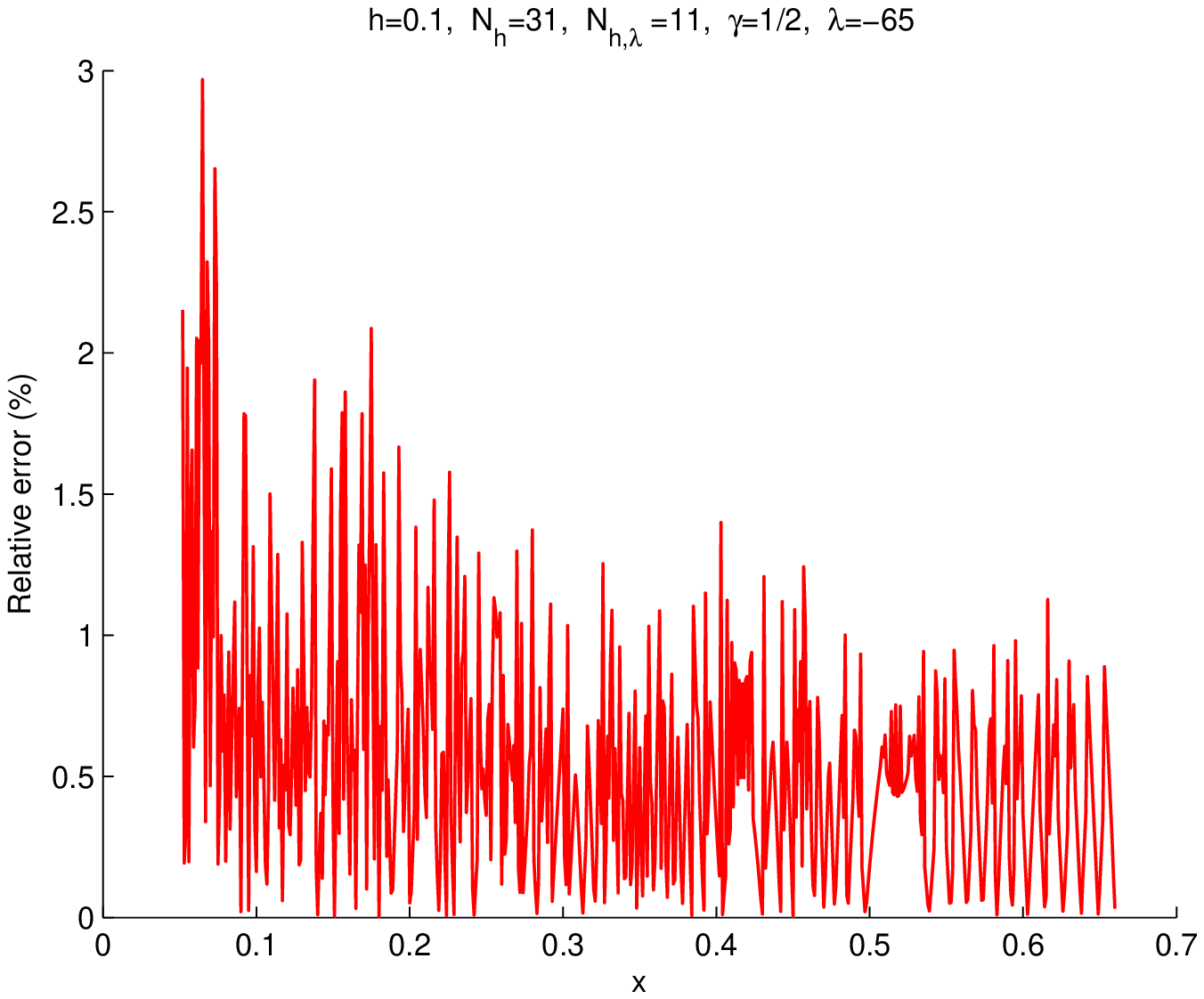}}
\caption{Reconstruction of a part of the ABP signal (a) and relative error (b) for $\lambda=-65$, $h=0.1$ and $\gamma=\frac{1}{2}$}
\label{pression-h01-lamda65-gamma05}
\end{center}
\end{figure}

\begin{figure}[htbp]
\begin{center}
\subfigure[]{\includegraphics[width=7cm]{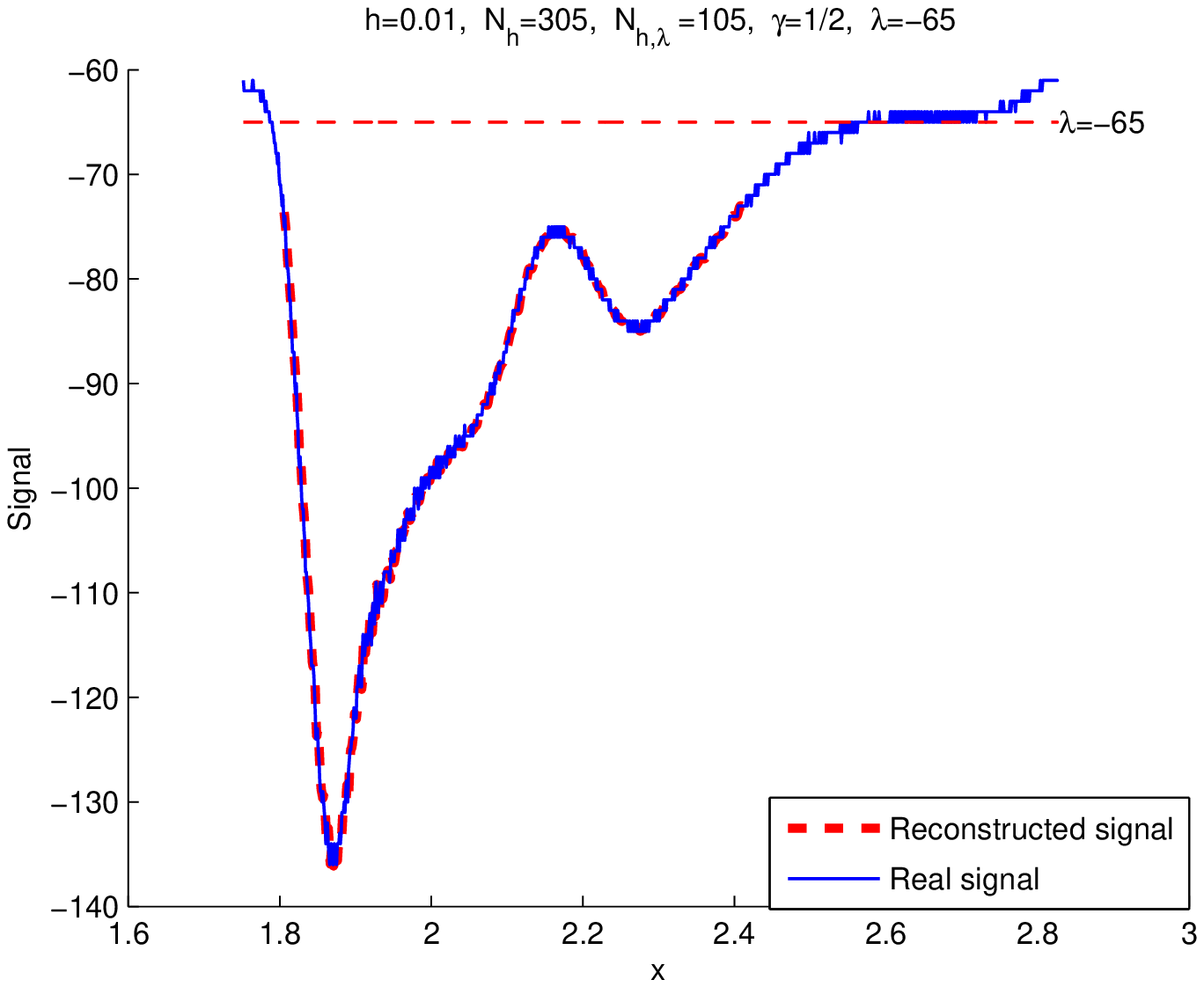}}
\subfigure[]{\includegraphics[width=7cm]{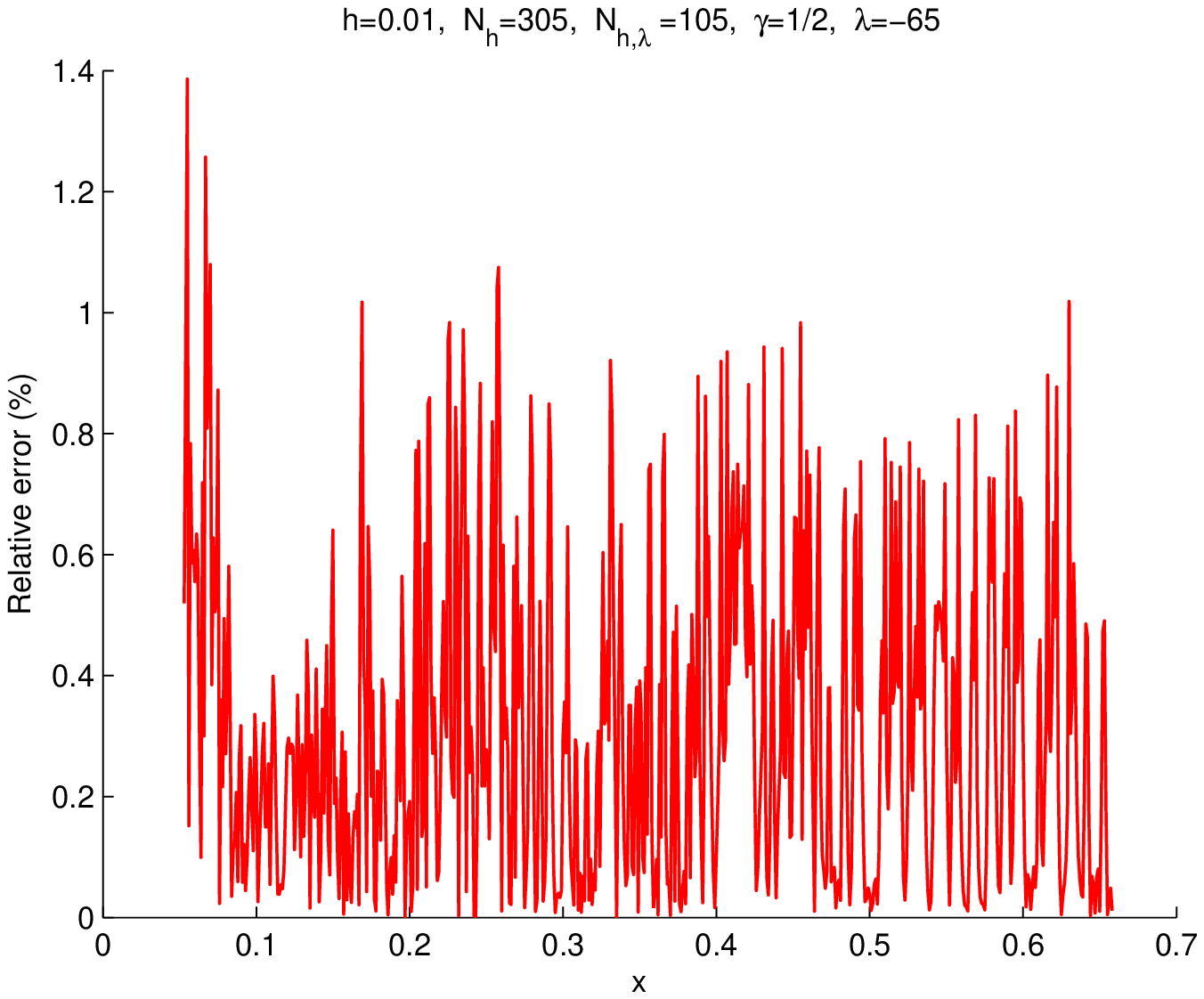}}
\caption{Reconstruction of a part of the ABP  signal (a) and relative error (b)  for $\lambda=-65$, $h=0.01$ and $\gamma=\frac{1}{2}$}
\label{pression-h001-lamda65-gamma05}
\end{center}
\end{figure}

\begin{figure}[htbp]
\begin{center}
\subfigure[]{\includegraphics[width=7cm]{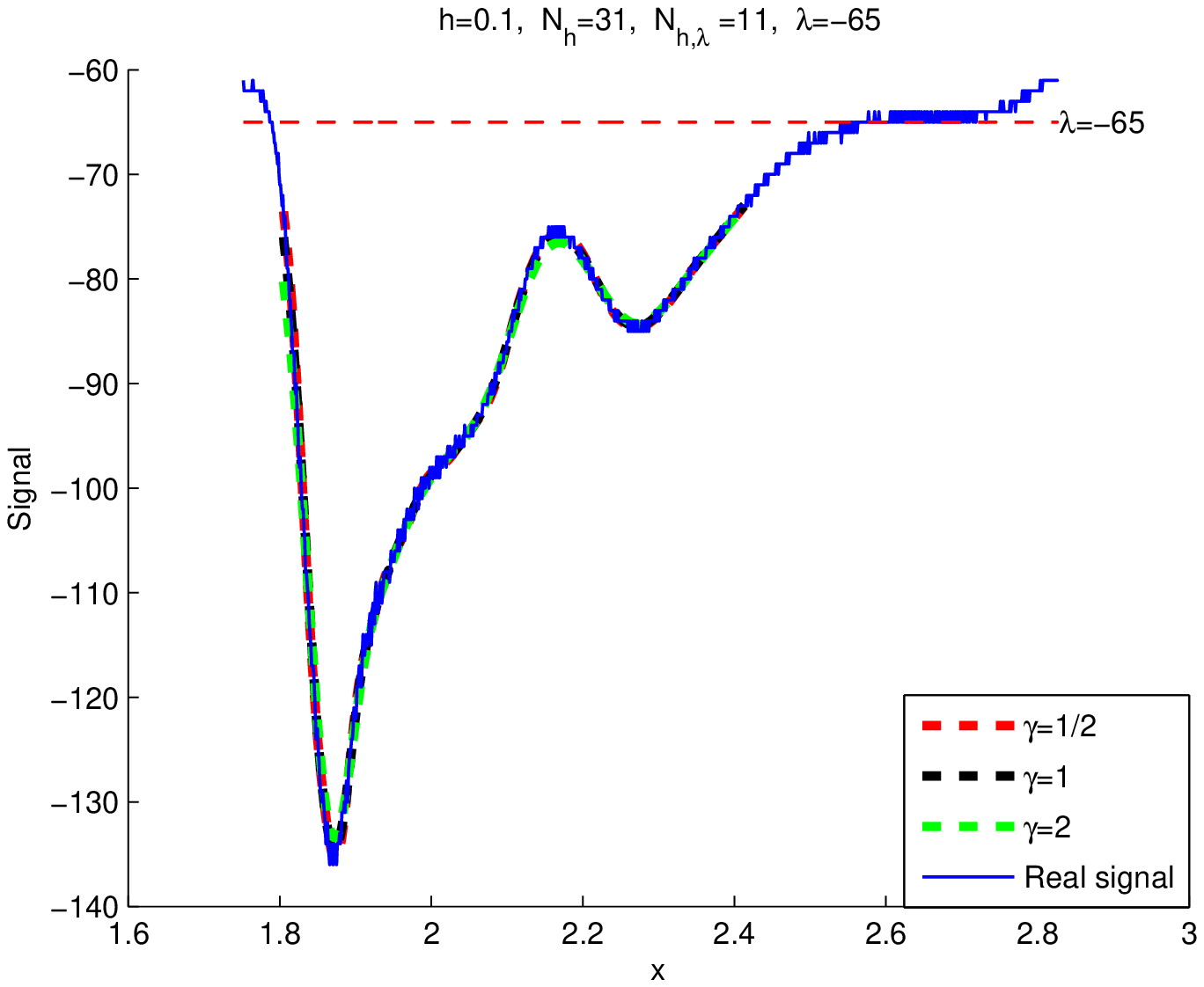}}
\subfigure[]{\includegraphics[width=7cm]{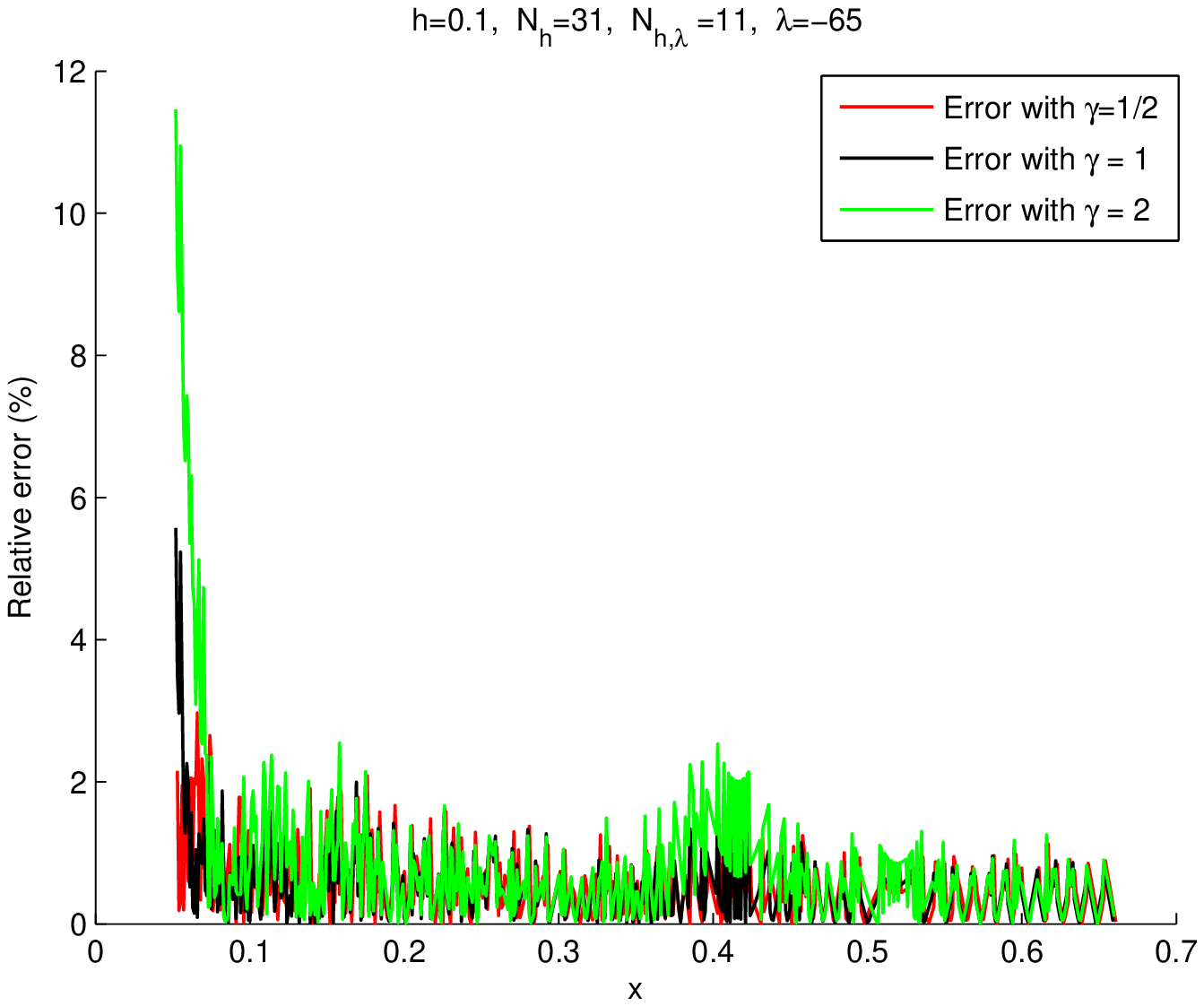}}
\caption{Reconstruction of a part of the ABP signal (a) and relative error (b) for $\lambda=-65$, $h=0.1$ and $\gamma=\frac{1}{2}, 1, 2$}
\label{pression-h01-lamda65-gamma05-1-2}
\end{center}
\end{figure}

\begin{figure}[htbp]
\begin{center}
\subfigure[]{\includegraphics[width=7cm]{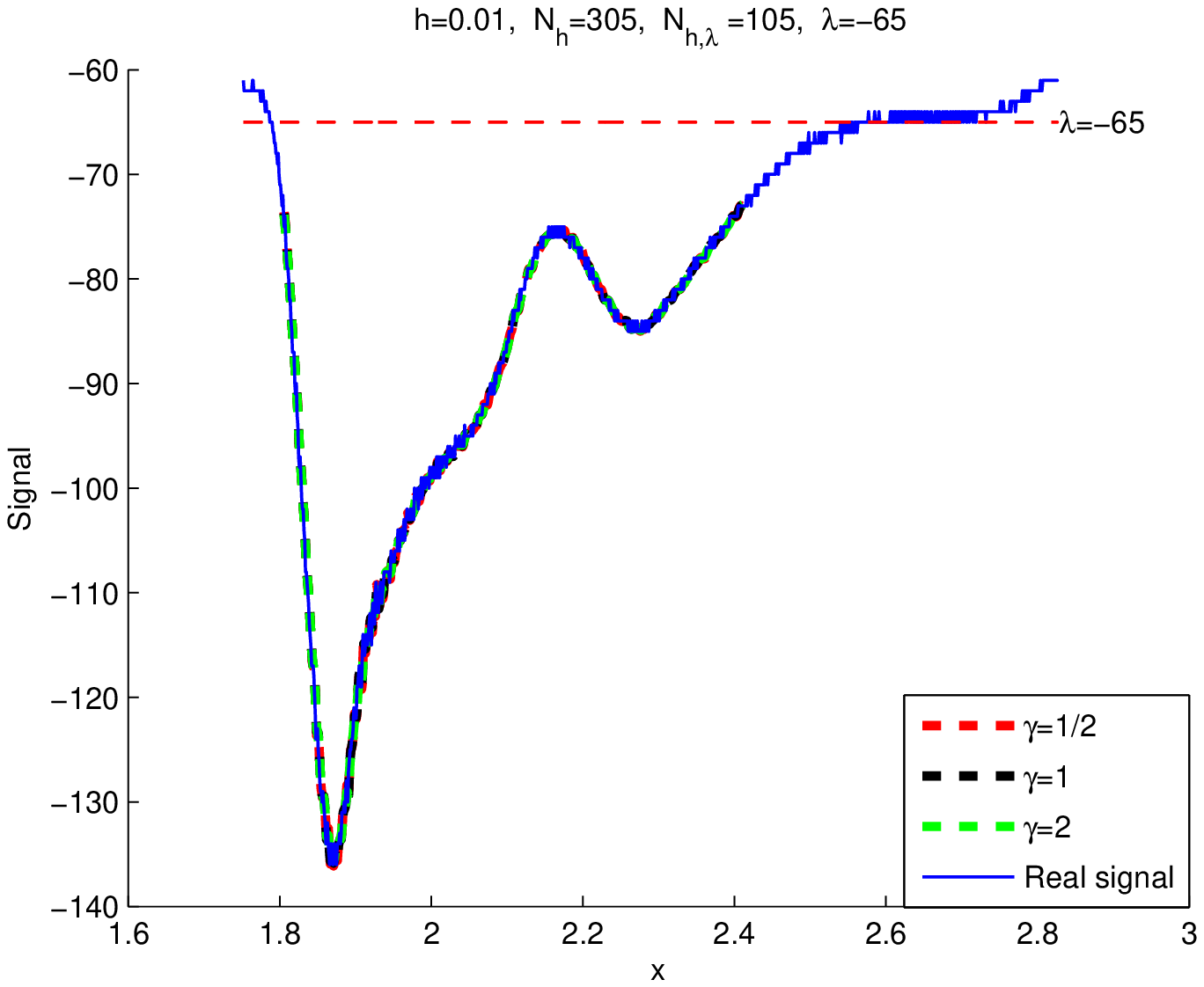}}
\subfigure[]{\includegraphics[width=7cm]{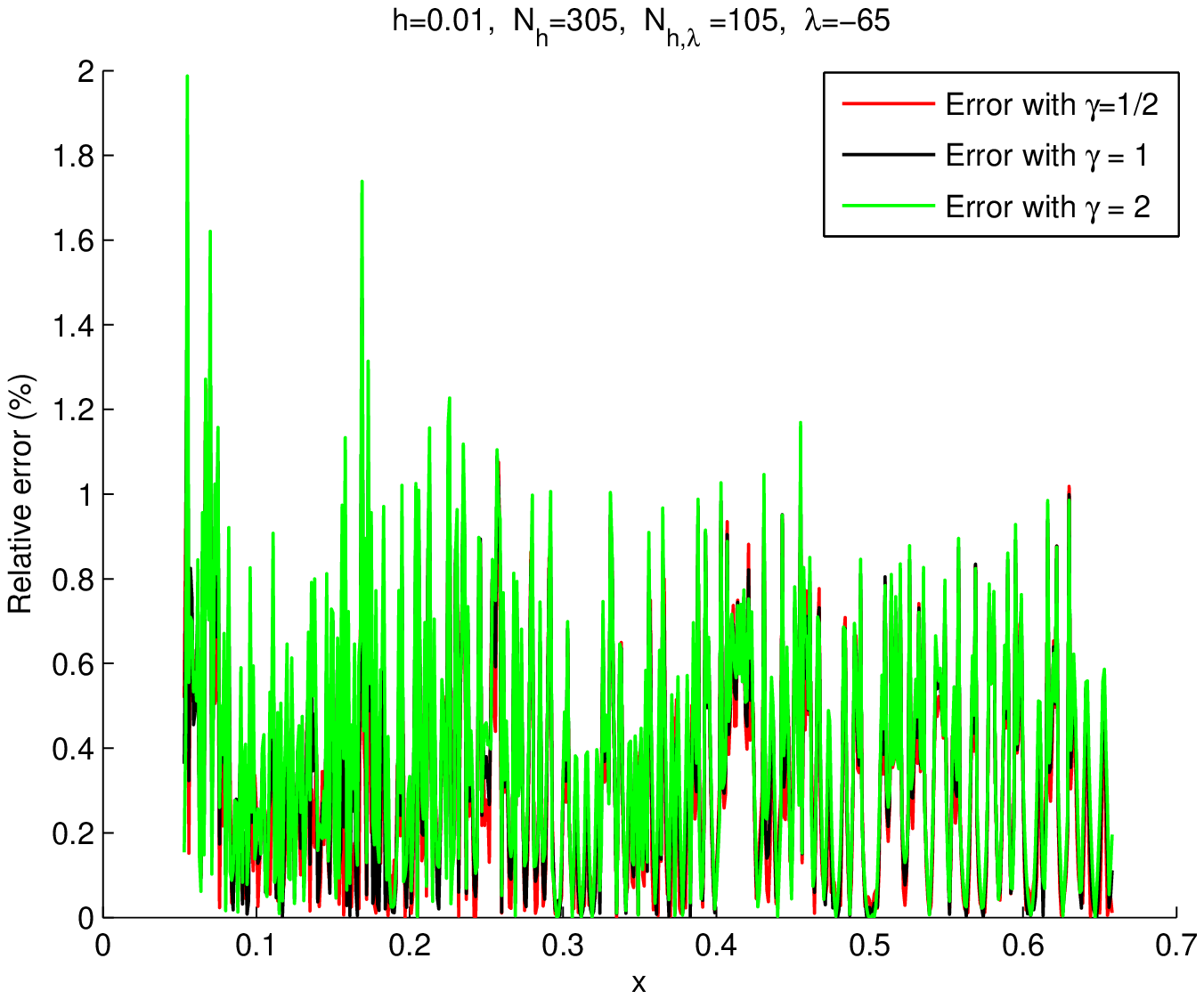}}
\caption{Reconstruction of a part of the ABP  signal (a) and relative error (b)  for $\lambda=-65$, $h=0.01$ and $\gamma=\frac{1}{2}, 1, 2$}
\label{pression-h001-lamda65-gamma05-1-2}
\end{center}
\end{figure}

Fig. \ref{pression-h01-lamda70-gamma05} and fig. \ref{pression-h001-lamda70-gamma05} as well as fig. \ref{pression-h01-lamda65-gamma05} and fig. \ref{pression-h001-lamda65-gamma05} clearly show that, as $h$ decreases, the approximation improves. Fig. \ref{pression-h01-lamda100-gamma05} and fig. \ref{pression-h01-lamda100-gamma05-1-2} suggest that, for  $h$ fixed, as  the distance of $\lambda$ to $- y(K)$ increases,  the estimate improves in $K$. Fig. \ref{pression-h01-lamda70-gamma05-1-2}, fig. \ref{pression-h001-lamda70-gamma05-1-2}, fig. \ref{pression-h01-lamda100-gamma05-1-2}, fig. \ref{pression-h01-lamda65-gamma05-1-2} and fig. \ref{pression-h001-lamda65-gamma05-1-2} illustrate the influence of $\gamma$. The differences between the errors of reconstruction in the three cases are not significant. However, the error for $\gamma=2$ is not regular as it was the case for the toy model. This  can be explained by the non regularity of our signal.   Note also that these reconstruction formulae appear as a regularization of the signal.

\begin{figure}[htbp]
\begin{center}
\subfigure[]{\includegraphics[width=7cm]{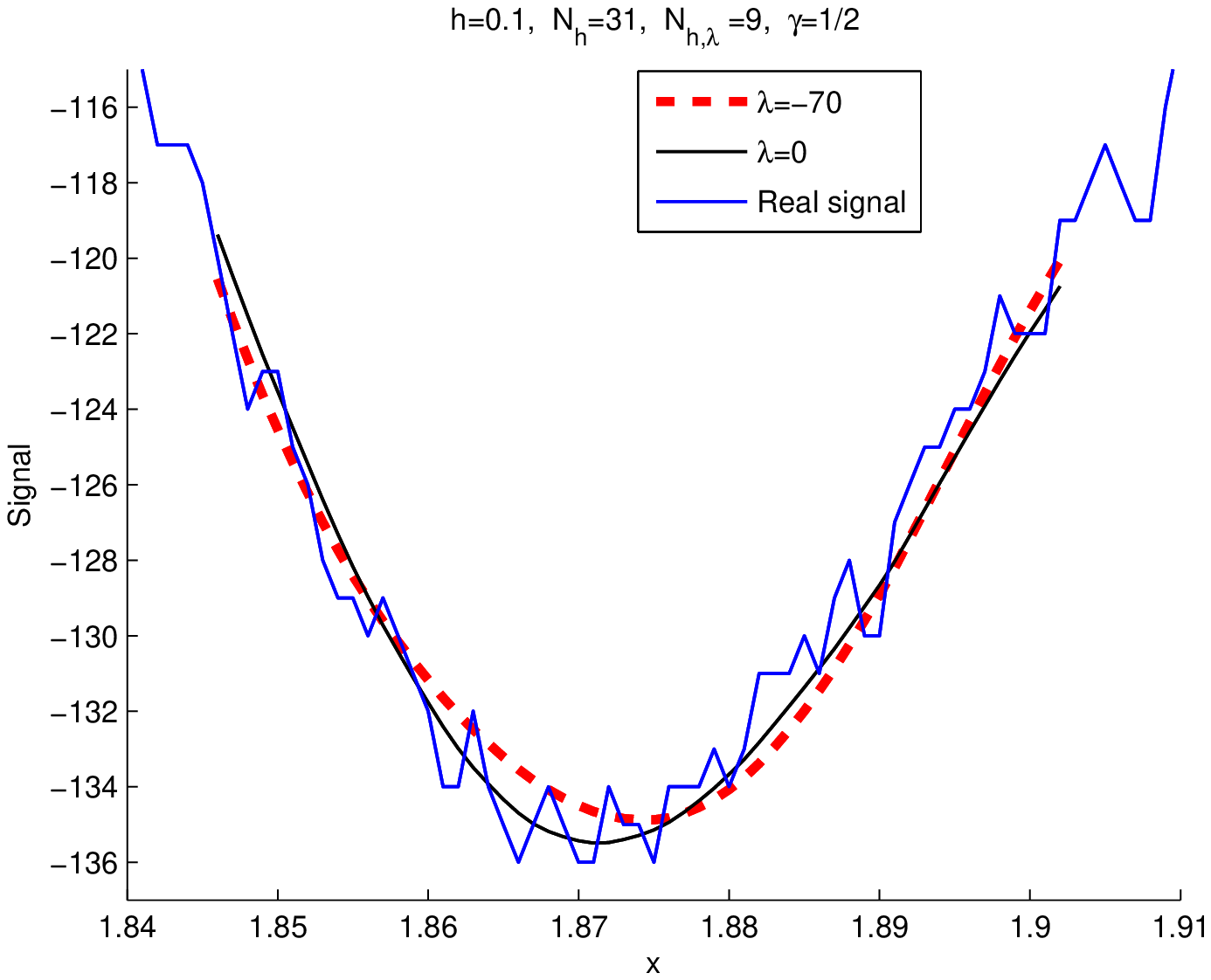}}
\subfigure[]{\includegraphics[width=7cm]{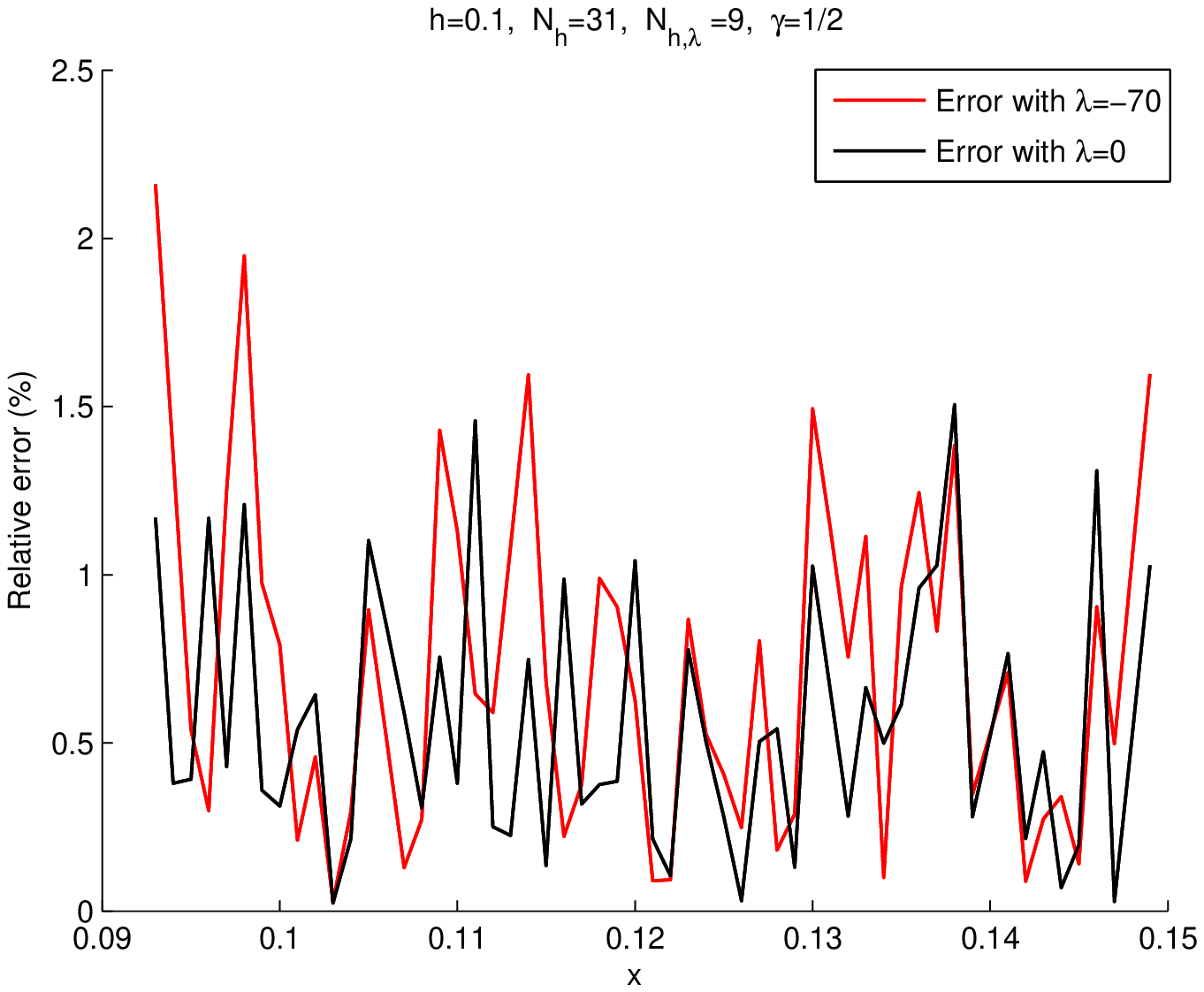}}
\caption{Reconstruction of a part of the ABP signal (a) and relative error (b) for $\lambda=0$ and $-70$, $h=0.1$ and $\gamma=\frac{1}{2}$}
\label{pression-h01-lamda0-70-gamma05}
\end{center}
\end{figure}

\begin{figure}[htbp]
\begin{center}
\subfigure[]{\includegraphics[width=7cm]{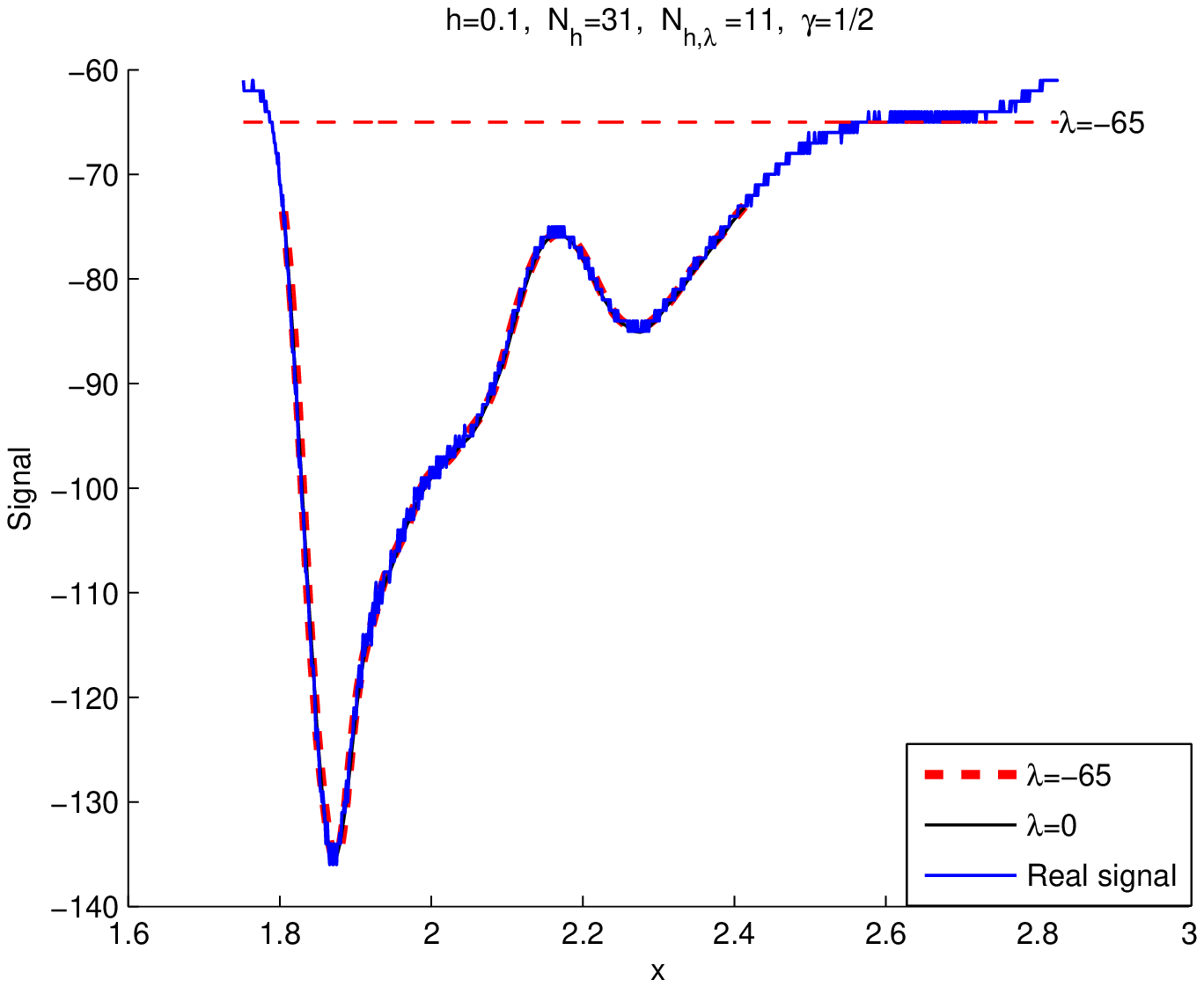}}
\subfigure[]{\includegraphics[width=7cm]{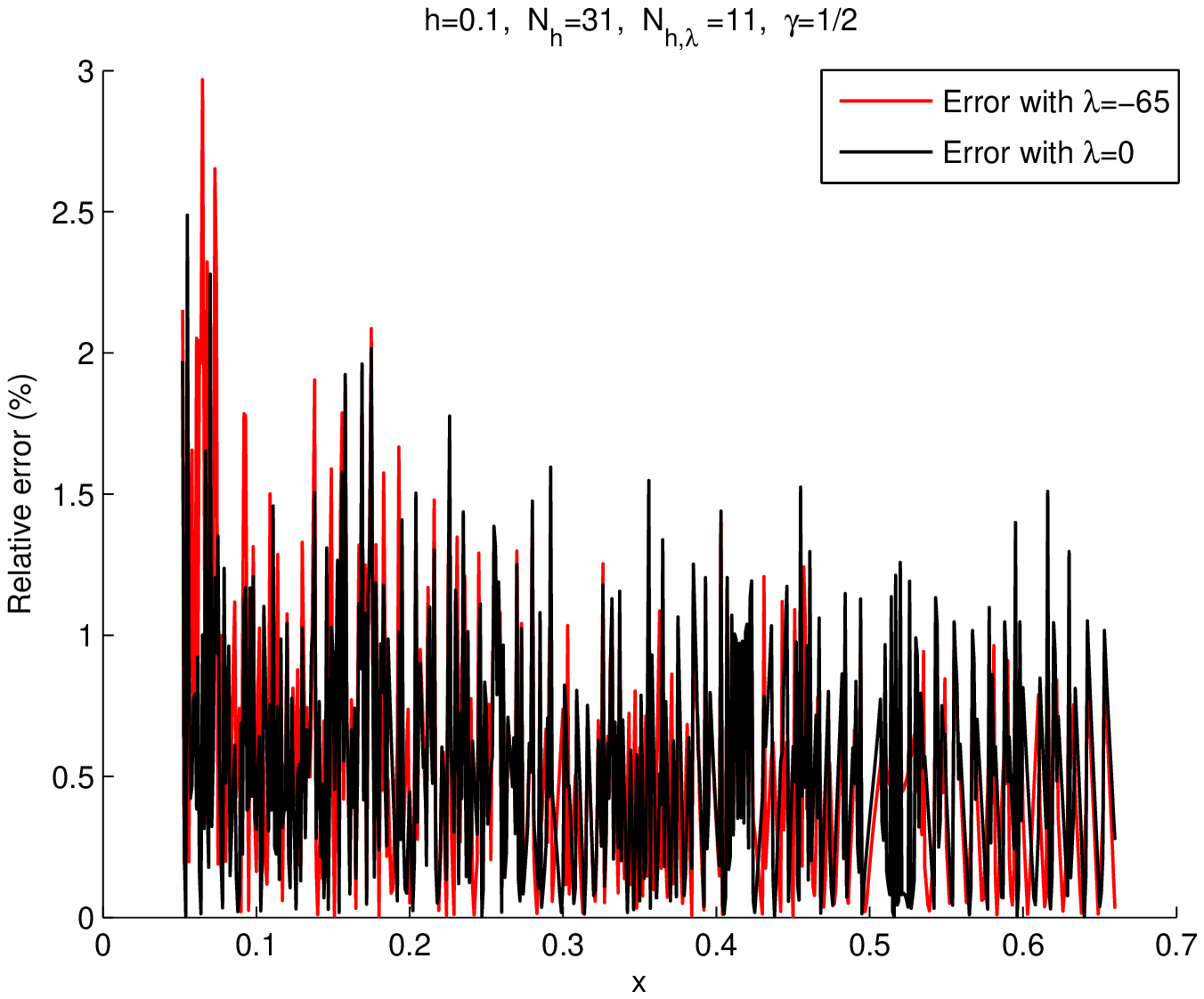}}
\caption{Reconstruction of a part of the ABP signal (a) and relative error (b) for $\lambda=0$ and $-65$, $h=0.1$ and $\gamma=\frac{1}{2}$}
\label{pression-h01-lamda0-65-gamma05}
\end{center}
\end{figure}

Fig. \ref{pression-h01-lamda0-70-gamma05}, fig. \ref{pression-h001-lamda0-70-gamma05}, fig. \ref{pression-h01-lamda0-65-gamma05} and fig. \ref{pression-h001-lamda0-65-gamma05} compare the case $\lambda =0$ to the case  $\lambda= -70 $ for the first example and to the case $\lambda =-65$ for the second example.  The errors are of the same order of magnitude. So it confirms our idea that we have not to take all the eigenvalues to estimate the signal.
Indeed,  we have found in the two cases a smaller  $\lambda$ such that the estimate in $K$ is satisfactory with a smaller number of eigenvalues and corresponding eigenfunctions to use.

\begin{figure}[htbp]
\begin{center}
 \includegraphics[width=10cm]{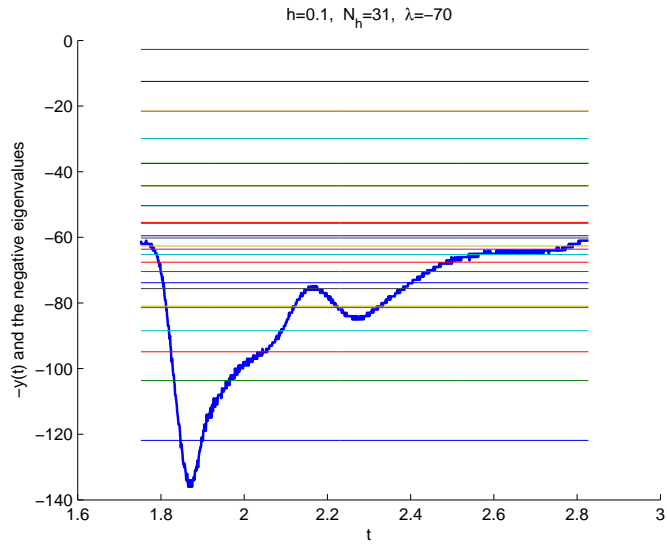}
 \caption{The negative eigenvalues of the Schr\"{o}dinger operator with a potential given by ABP signal }
 \label{VAP_lamda70}
\end{center}
\end{figure}

Fig. \ref{VAP_lamda70}  is a good illustration of many semi-classical properties. First we see that, below the energy $-60$,  there is a concentration of the eigenvalues near the critical points of $-y$.  Far above $-60$, the behavior is asymptotically given (for fixed $h$) by the eigenvalues of
the  periodic realization of $-h^2\dfrac{d^2}{d x^2}$ in $]a,b[$ with $a$ and $b$ describing the beginning and the end of an arterial blood pressure beat respectively.

\section{Conclusion}

In continuation of \cite{Laleg:08}, \cite{LaCrSo:10}, \cite{LaMePaCoVa:10}, we have explored the possibility of reconstructing  the signal $y$ using
 spectral quantities  associated with some self-adjoint realization of an $h$-dependent Schr\"odinger operator ($h>0$) $-h^2\frac{d^2}{dx^2} -y (x)$, the parameter $h$  tending to $0$. Using on one hand theoretical results in semi-classical analysis and on the other hand numerical computations, we can formulate
  the following remarks.
 \begin{itemize}
 \item As $h$ ($h>0$) decreases (semi-classical regime), the approximation of  $y$ by $y_{h,\gamma}(x,\lambda)$  improves.
 \item Semi-classical analysis suggests also to take $\gamma \geq 1$.
 \item However, the number of negative eigenvalues to be computed increases as $h$ becomes smaller. So the numerical computations become difficult if $h$ is too small.
 \item For a given interval $K$, a clever choice of $\lambda$ makes  possible to get a good approximation of the signal in $K$
    in the semi-classical regime with a smaller number of eigenvalues.
    \item  Finally, the numerical computations suggest also (but we are outside the theoretical considerations of our paper)  that one can consider
    a choice of  $h$ (not necessarily small)   and hope a good reconstruction  of the signal by choosing  appropriate values of $\lambda$ and $\gamma$. This should be the object of another work.
    \end{itemize}

\paragraph{Acknowledgements:}  The authors would like to thank G. Karadzhov for useful discussions around his work and Doctor Yves Papelier from Hospital Béclère in Clamart for providing us arterial blood pressure data.

\end{document}